%% file: Novelrobin.tex
\documentclass[a4paper,12pt,reqno]{amsart}
\usepackage{amssymb}
\usepackage{delarray}

\newtheorem{thm}{Theorem}[section]

\newtheorem{prop}[thm]{Proposition}
\theoremstyle{definition}
\newtheorem{defn}[thm]{Definition}
\theoremstyle{remark}
\newtheorem{rem}[thm]{Remark}
\numberwithin{equation}{section}
\newtheorem{example}{Example}[section]

\newcommand{\beq}  {\begin{equation}}
\newcommand{\eeq}  {\end{equation}  }

\newcommand{\bit}  {\begin{itemize}}
\newcommand{\eit}  {\end{itemize}  }

\newcommand{\ben}  {\begin{enumerate}}
\newcommand{\een}  {\end{enumerate}  }
\newcommand{\eqnref}[1]{(\ref{#1})}
\newcommand{\bpr}  {\begin{proof}}
\newcommand{\epr}  {\end{proof}  }

\newcommand{\Real}{\mathbb R}

\def\l{\lambda}

\def\R{{\mathbb R}}

\newcommand{\pint} {\mathrm{P}\int_{0}^\infty}
\begin{document}

\title[Novel Exact Solutions]{Novel Exact Solutions for PDEs with Mixed Boundary Conditions }%

\author{Mark Craddock, Martino Grasselli and  Andrea Mazzoran}

\address{Author affiliations: 
\newline Craddock: School of Mathematical and Physical Sciences, University of
Technology Sydney, PO Box 123, Broadway, New South Wales 2007,
Australia. Email: Mark.Craddock@uts.edu.au. 
\newline Grasselli, corresponding author: Department of Mathematics ``Tullio Levi Civita'',
University of Padova, via Trieste 63, 35121 Padova, Italy, and L\'eonard de Vinci P\^ole Universitaire, Research Center, 92916 Paris la D\'efense. Email: grassell@math.unipd.it. 
\newline Mazzoran: Department of Mathematics ``Tullio Levi Civita'',
University of Padova, via Trieste 63, 35121 Padova, Italy. Email: andreamazzoran@hotmail.it}

\date{November 16, 2023}%
\begin{abstract}
We develop methods for the solution of inhomogeneous Robin type boundary value problems (BVPs) that arise for certain linear parabolic Partial Differential Equations (PDEs) on a half line, as well as a second order generalisation. We are able to obtain non-standard solutions to equations arising in a range of areas, including mathematical finance, stochastic analysis, hyperbolic geometry and mathematical physics. Our approach  uses the odd and even Hilbert transforms. The solutions we obtain and the method itself seem to be new.
\end{abstract}
\maketitle
\input  Robinproblems
\bibliographystyle{plain}
\bibliography{robbiblio}
\end{document}

%% file: Robinproblems.tex
\section{Introduction}

There is a well established theory for the solution of parabolic PDEs subject to the most common types of boundary conditions. The book by Friedman \cite{Fried64}, provides a rigorous introduction to this topic. The classical method makes use of a fundamental solution of the PDE which also satisfies the boundary conditions. The construction is straightforward and we will present an illustrative example in Appendix 1.

However the classical method can produce cumbersome representations of the solution and the desired fundamental solution may not even be known. So we ask if it is possible to  construct an analytical solution to a boundary value problem, using only   elementary solutions, without using the fundamental solution required by the classical theory? This has potentially important practical implications, because there are many PDEs for which elementary solutions are readily obtainable, but for which appropriate fundamental solutions are not known.

We will focus on parabolic PDEs on a half line $(b,\infty)$, subject to the boundary condition
\begin{equation}\label{genrobcond}
\alpha u(b,t)+\beta u_x(b,t)+\gamma u_{xx}(b,t)=g(t).
\end{equation}
We will refer to this as a second order Robin condition. If $\gamma=0$ this reduces to the usual Robin condition.\footnote{There is an historical curiosity here. The boundary condition with $\gamma=0$ is named after Gustave Robin. However according to \cite{GA98}, Robin never stated or studied this type of boundary condition. It appears nowhere in his collected works. Why the condition came to be named after Robin is apparently a mystery.} The Dirichlet condition $u(b,t)=0$ and the Neumann condition $u_x(b,t)=0$ are also special cases. We are also able to solve certain moving boundary problems. We will give an example in the final section.

We develop a new method for the solution of these problems. Our technique relies on the odd and even Hilbert transforms and does not require a fundamental solution. We only need two elementary solutions. These can be found in a number of ways, such as separation of variables.

Suppose that we have a linear parabolic PDE $u_t=Lu$ on some interval $(b,\infty) $. $L$ may be time dependent. The essential idea is to look for solutions of the form
\begin{align}\label{formsoln}
u(x,t)=\int_0^\infty\varphi(\xi)w_1(x,t;\xi)d\xi+\int_0^\infty\psi(\xi)w_2(x,t;\xi)d\xi,
\end{align}
where $w_{1,2}(x,t;\xi)$ satisfy the PDE for each $\xi$. Imposing \eqref{genrobcond} and an initial condition leads to a pair of integral equations for $\varphi$ and $\psi.$

In general there is no reason to suppose that these equations will be analytically tractable, though one might attempt to solve them numerically. This idea is  reminiscent of the well known boundary integral method used for higher dimensional BVPs. See for example \cite{MCos90}.

However it turns out that for certain types of important problems, these equations admit explicit solutions and lead to  a representation of the solution that differs  from those produced by the classical method. In particular, we have a solution of the BVP which does not rely upon us knowing a fundamental solution.

In the current work we focus on problems where $w_1, w_2$ involve the sine and cosine functions. This will be made explicit in Section \ref{DaveHilb}. We will develop this new  technique and present an interesting selection of examples. The outline of the paper is as follows:
%

We begin with a general discussion of boundary value problems and the representation of their solutions. We mention some recent work, particularly of Fokas.

Following this, we turn to the solution of the problem that motivated this study, namely the second order Robin problem for the Black-Scholes equation.
This is the content of Section \ref{SecordRobBS}.
We reduce the solution of this problem to the inversion of a Fourier sine transform. Our result appears to be new. We also briefly give the fundamental solution for the classical Robin problem, which can be solved by the same method.

Since the fundamental solution in the second order case is cumbersome, we ask if it is possible to solve the Robin problem without the fundamental solution?
The answer is yes and we develop our new theory for the explicit case of the Black-Scholes equation. The method relies on the odd and even Hilbert transforms and we begin by presenting the properties of these transforms that we need. See \ref{oeht}. This work begins in Section \ref{DaveHilb} with the classical Robin problem. We then  proceed to the second order Robin problem in Section \ref{DaveHilb2}. We follow this with an example.

Then we turn to the second order Robin problem for a larger class of PDEs. This is in Section \ref{QE2d}. Theorem \ref{AGentheorem} gives the explicit solution of a class of second order Robin problems with solutions involving the sine and cosine problems. The techniques we developed for the Black-Scholes example make the proof easier. This new method is the main contribution of our paper.

A set of examples of PDEs and the families of elementary solutions which our method requires are then presented. This is Subsection \ref{exes}. Following this we turn to the question of PDEs which are not covered by Theorem \ref{AGentheorem}. We briefly mention PDEs with time dependent coefficients then solve two further problems. First for the harmonic oscillator, which requires a different set of elementary solutions. This is in Subsection \ref{LT}. In Section \ref{BESQ5} we solve the Robin problem for a five dimensional squared Bessel process. This has the feature of using elementary solutions that involve sums of sines and cosines. In both cases our method is effective.

We conclude with a brief discussion of open problems and future directions.

\subsubsection{Boundary Value Problems}

The literature on boundary value problems is enormous, but  Sagan's book provides a good elementary introduction \cite{Sag61}. Boundary integral methods play an important role in the study of multi-dimensional BVPs.  See for example McLean's book \cite{McLean2000}. For a study of boundary conditions in the theory of diffusions, see \cite{SV72}. We also mention the work of Fokas (e.g. \cite{Fok}). This presents a novel integral transform method for the solution of BVPs on a half line. There are obviously thousands of references and we could not attempt an exhaustive list.

However, for Robin problems we mention some recent work. In Abels and Moser \cite{abels2022convergence},  a nonlinear Robin boundary condition in a bounded smooth domain is investigated.  Geng and Zhuge \cite{geng2020oscillatory}, study a family of second-order elliptic systems subject to a periodically oscillating Robin boundary condition. In \cite{lawley2015newa},  the authors study the diffusion equation with a stochastic boundary which randomly switches between a Dirichlet and a Neumann condition, proving that the mean of the solution satisfies a new type of Robin condition.

It is worth mentioning also the paper \cite{daners2008uniqueness}, where the authors consider the eigenvalues of the Robin boundary value
problem for the Laplacian. Finally Bondurant and Fulling \cite{BF2005}, introduce a map between Dirichlet and Robin boundary conditions for linear constant coefficient equations.

The current work arose from the study of barrier options within the Black-Scholes (BS) framework. For so called {\it knock out barrier options}, we  have absorbing boundaries. Reflecting boundaries occur where the option is {\it knocked-in}. Although many formulae exist for the pricing of barrier options, the solution of the Robin problem for the BS equation does not seem to have appeared in the literature. We actually solve the second order Robin problem for the BS equation.

The fundamental solution that we obtain appears to be new.  A glance at Appendix 1 shows that this  solution is extremely complex. This leads to the question which is the main part of our study: Is it possible to solve the BVP without the fundamental solution?  Although  any well posed BVP has a unique solution, there is in general no unique {\it representation}  of that solution. This is not a trivial fact. For example, consider the BVP
\begin{align*}
u_t&=u_{xx},\     x> 0,t>0, u(x,0) =f(x),x>0, f\in L^1(\Real^+),\\
 u(0,t)&-\gamma u_x(0,t)=g(t),\ t>0.
\end{align*}
The solution of this can be written in terms of the classical heat kernel. This is well known. See Cannon's book \cite{Cannon84} for the details.

However another representation of the solution to this problem   was obtained by Fokas. We quote the result from \cite{Fok}. The solution can be written
\begin{align*}
u(x,t)&=\frac{1}{2\pi}\int_{-\infty}^\infty e^{ikx-k^2t}\widehat{f}(k)dk-\frac{1}{2\pi}\int_{\partial D^+}e^{ikx-k^2t}\left[\frac{2k}{k+i\gamma}G_R(k^2)\right. \\&\left. -k\frac{k-i\gamma}{k+i\gamma}\widehat{f}(-k)\right]dk+\Theta(-\gamma)2\gamma e^{\gamma x+\gamma^2t}[G_R(-\gamma^2)-\widehat{f}(i\gamma)].
\end{align*}
Here $\Theta$ is the Heaviside step function, $G_R(k)=\int_0^T e^{ks}g(s)ds,\ k\in\mathbb{C}$, $D^+$ is the wedge in the upper half plane making an angle of $\pi/4$ with the real axis on both sides of the origin (see figure 2, p4 of \cite{Fok}) and $\widehat{f}(k)=\int_0^\infty f(x)e^{-ikx}dx.$

A discussion of the advantages of this alternative way of obtaining the solution is beyond our scope, though we remark that Fokas type representations can often be evaluated numerically with considerable efficiency. However this is a large subject and for brevity we refer the reader to the aforementioned book by Fokas  and the works of Donaldson (e.g. \cite{JADon75}). There is a considerable literature on this topic.

%

We will use a combination of separation of variables and classical transform methods.  The Hilbert transform has been used to solve integral equations for over a hundred years. See \cite{Kin2009b} for a lengthy discussion with examples. We actually  use  the related {\it odd} and {\it even} Hilbert transforms. Our method has the elegant feature that it turns a pair of integral equations into a pair of simultaneous equations.

Some of our results rely upon the inversion of a Laplace transform. However there are thousands of Laplace transform pairs known
and more can be constructed by standard methods, see \cite{RK66}. For polynomial data, inversion produces Dirac delta functions and their derivatives and these are easy to handle. If $g(t)=0$, then there is no Laplace transform to invert.
\begin{rem}
We make an important comment here. There are many PDEs on the line that can be mapped to the heat equation. Some of our examples have this property, though not all. (For example equation  \eqref{5dsbes} cannot be reduced to the heat equation).  One might think that it would be more efficient to reduce a suitable PDE to the heat equation, solve the resulting BVP, then map back. However this ignores the question of what happens to the boundary conditions under the change of variables? This is a crucial question. The resulting boundary value problem may not have a known solution.

To illustrate, suppose that we wish to solve $u_t=u_{xx}+Axu$, $u(x,0)=f(x)$, $u(0,t)=0$. The PDE can be mapped to the heat equation. However the result of the mapping is a formidable {\it moving boundary problem}. The complete details are in \cite{CG19}. Thus mapping to the heat equation in this case makes the problem much harder.

A similar phenomenon occurs with the harmonic oscillator. Reducing even the simple problem $u_t=u_{xx}-x^2u$, $u(x,0)=f(x)$, $u(0,t)=0$ to the heat equation produces another moving boundary problem, which is much harder than the original problem. There are numerous examples of this phenomenon.

There are also many PDEs which can theoretically be mapped to the heat equation, but the change of variables is itself impossible to compute. Consider the equation $u_t=\sigma(x)u_{xx}+f(x)u_x, x\in\Omega\subset\mathbb{R}.$ Suppose that it can be reduced to the heat equation. To do so we first let $y=\int^x_{x_0}(\sigma(z))^{-1/2}dz$. This makes the coefficient of the second derivative term equal to 1.  For arbitrary $\sigma$, there is no reason why this integral should be computable. Inverting the change of variables to write $x$ in terms of $y$ may also be impossible.

Thus methods which produce solutions of BVPs without the need for a change of variables are important. Our results  yield  novel solutions for a wide variety of problems without having to make a change of variables.

\end{rem}

\section{The Second Order Robin Problem for the Black-Scholes Equation}\label{SecordRobBS}
For the theory of option pricing and stochastic calculus, we refer the  reader to a standard reference such as \cite{DHW95}.
The conventional method for studying the Black-Scholes equation is to reduce to the heat equation,
but we work in the original variables. For our purposes the
second order Robin problem can be written

\begin{align}\label{Robin Problem}
    w_t &= \frac{1}{2}\sigma^2S^2 w_{SS}+rSw_S,\\
    w(S,0)&=f(S), \label{Robin Problemb} \\
\alpha w(b,t)+\beta w_S(b,t) + \gamma w_{SS}(b,t) & =g(t)  \label{Robin Problemc},
\end{align}
where $S>b>0$.

We will assume that the solution is nonnegative and satisfies a bound of the form $w(S,t)\leq MS^\theta$ for some positive constants $M,\theta.$ This is the usual type of behaviour as $S\to\infty$ that appears in the literature. See for example the discussion of boundary conditions for the Black-Scholes equation in \cite{DHW95}.

We will construct a fundamental solution for this problem in the case $g(t)=0$.
We first obtain separable solutions which satisfy the boundary condition. We then use these to construct a solution
which also satisfies the initial condition. This second problem can be reduced to the solution of a tractable integral equation.

We make the ansatz $w(S,t)=e^{\lambda t}v(S)$. Then $\frac{1}{2}\sigma^2S^2 v''(S)+rSv'(S)-\lambda v(S)=0, $
with  the constant  $\lambda$ and  the function  $v$ to be determined.
Assuming $v(S)=S^\delta, $
we obtain the condition.
\beq
    \frac{1}{2}\sigma^2S^2 \delta (\delta-1)S^{\delta-2}+rS\delta S^{\delta -1}-\lambda S^\delta=0.
\nonumber\eeq
This means that we must have $\frac{1}{2}\sigma^2 \delta (\delta-1)+r\delta -\lambda =0,$
which gives

\beq
    \delta=\frac{-\left(r-\frac{\sigma^2}{2}\right)\pm \sqrt{\left(r-\frac{\sigma^2}{2}\right)^2+2\sigma^2\lambda}}{\sigma^2}.
\nonumber\eeq

We now set $\left(r-\frac{\sigma^2}{2}\right)^2+2\sigma^2\lambda=-\xi^2 \sigma^4.$
So we have
\begin{equation}
\label{lambda definition}
    \lambda= -\frac{\xi^2 \sigma^2}{2}-\frac{1}{2\sigma^2}\left(r-\frac{\sigma^2}{2}\right)^2.
\end{equation}
Hence we can write $\delta =\mu\pm i\xi,$ where $\mu=\frac{\sigma^2-2r}{2\sigma^2}.$

This gives us the solution which one can easily check is valid for all $\xi>0.$ This is
\begin{align}
    w(S,t,\lambda)&=e^{\lambda t}\left( c_1  S^{\mu+i\xi}+c_2 S^{\mu-i\xi}\right).
\nonumber\end{align}
We now wish to construct a solution of the PDE that satisfies the  boundary condition. Observe that
\begin{align*}
    S^{\mu\pm i\xi}=S^{\mu}\left( \cos (\xi \ln S )\pm i \sin (\xi \ln S )\right).
\end{align*}
The real and imaginary parts must both satisfy the equation. So we obtain a  solution that can be written as
\begin{equation}
\label{PDE general solution}
    \begin{aligned}
        w^\xi (S,t)&= e^{\lambda t}S^{\mu}\left( A\cos (\xi \ln S )+B \sin (\xi \ln S )\right)\\        &=e^{-\frac{\xi^2 \sigma^2}{2}t+ct}S^{\mu}\left( A\cos (\xi \ln S )+B \sin (\xi \ln S )\right),
    \end{aligned}
\end{equation}
where $c=-\frac{1}{2\sigma^2}\left(r-\frac{\sigma^2}{2}\right)^2$ and $A,B$ are constants. We choose $A$ and $B$ by requiring that the solution satisfies the second order Robin boundary condition  with $g(t)=0$.
Substituting the solution into the PDE and imposing the boundary condition, we obtain 

\begin{align}
A&=\sin (\xi  \log (b))z_1+\xi z_2 \cos (\xi \log (b)), \label{generalA}\\
B&=-\cos (\xi  \log (b))z_1+\xi z_2\sin (\xi  \log (b)),\label{generalB}
\end{align}
where $z_1= \left(\alpha  b^2+\beta  b \mu +\gamma  (\mu -1) \mu -\gamma  \xi ^2\right)$ and $z_2= (b \beta +\gamma  (2 \mu -1)).$

The reader can check that for this choice of the coefficients $A$ and $B$,  the function $w^\xi$    solves the PDE  and it also satisfies the homogeneous form of the second order Robin boundary condition \eqref{Robin Problem}. These choices are not unique. However every other choice actually leads to the same solution. This is not hard to check, but it is somewhat tedious.

We now have to obtain a solution that also satisfies the initial data as well as the boundary condition.
To do this we will construct a solution of the form
\begin{align}
\label{General solution of PDE and Robin condition}
    w(S,t)=\int_0^\infty \varphi (\xi) w^\xi (S,t)d\xi.
\end{align}
If the function $\varphi$ has sufficient decay, it is easy to show that $w(S,t)$ is also a solution of our PDE. See \cite{Cra2008} for more on this idea.
We note that $w^{\xi}(S,t)$ is locally integrable in $\xi$, since it is continuous in $\xi$ for all $S>0$ and has Gaussian decay for all $t>0.$

Moreover, by construction,  \eqref{General solution of PDE and Robin condition}  satisfies the boundary condition. Our task now is to choose the function $\varphi$ so that the initial condition is also satisfied.
Taking $t=0$, we get the following integral equation
\begin{align}
    f(S)=\int_0^\infty \varphi (\xi) w^\xi (S,0)d\xi.
\nonumber\end{align}
That is,
\begin{align}
    \label{Integral Equation}
    \int_0^\infty \varphi (\xi)\left( A\cos (\xi \ln S )+B \sin (\xi \ln S )\right)
    d\xi=S^{-\mu}f(S).
\end{align}
This can be reduced to a Fourier sine transform. To see this, notice that using some elementary trigonometric identities, we can rewrite \eqref{Integral Equation} as
    \begin{align}
        \int_0^\infty \varphi (\xi)\left( \tilde{\beta} \xi \cos \left(\xi \ln \left(\frac{S}{b} \right)\right) + \left(\gamma \xi^2 - \tilde{\alpha}\right) \sin \left(\xi \ln \left( \frac{S}{b} \right)\right) \right)
        d\xi=F(S),
    \nonumber\end{align}
where  we set $F(S)=S^{-\mu}f(S) $  and
\begin{align}
	\tilde{\alpha}&=\alpha b^2 + \beta b \mu +\gamma \mu (\mu -1), \label{a tilde}\\
	\tilde{\beta}&=b \beta + \gamma (2 \mu -1).\label{b tilde}
\end{align}
We now reduce  this integral equation to a Fourier sine transform via the solution of a second order
ordinary differential equation.

Let $\zeta(x)=\int_0^\infty \varphi (\xi) \sin (\xi x)d\xi$. After making the change of variable $x=\ln \left( \frac{S}{b} \right)$ and differentiating twice, we see that $\zeta$ must satisfy the equation
\begin{equation}
\label{ODE for Boundary Condition}
	-\gamma \zeta^{\prime \prime}(x) + \tilde{\beta} \zeta^{\prime}(x) - \tilde{\alpha} \zeta(x) = b^{-\mu} e^{- \mu x}f(b e^x),
\end{equation}
with initial conditions
$\zeta(0)=0$ and $\zeta''(0)=0$. We will assume that  $\Delta= \tilde{\beta }^2-4 \gamma  \tilde{\alpha }>0$. For the special cases where $\Delta\leq 0$, the differential equation for $\zeta$ has different solutions. However we can proceed as we do here and we obtain fundamental solutions valid for those particular choices of parameter. We will omit the details for brevity.

Using variation of parameters, we see that the solution of \eqref{ODE for Boundary Condition} is given by

\begin{align*}
\zeta(x)=& \frac{2b^{-\mu }}{\tilde{\beta } \sqrt{\Delta}}  \left( \gamma f(b) e^{\frac{\tilde{\beta }x}{2 \gamma }} \sinh \left(\frac{x \sqrt{\Delta}}{2 \gamma }\right) - \tilde{\beta} \int_0^x f\left(be^z\right) \mathcal{K}(x,z)dz \right),
\end{align*}
where we set $\mathcal{K}(x,z)=e^{- \mu z+(x-z)\frac{\tilde{\beta }}{2 \gamma }}
\sinh \left(\frac{(x-z) \sqrt{\Delta}}{2 \gamma }\right).$

Recall that the inverse Fourier sine transform of $\widehat{f}\in L^1([0,\infty))$ is
\begin{align}
f(x)=\frac{2}{\pi}\int_0^\infty \widehat{f}(y)\sin(yx)dy,
\end{align}
 see \cite{Z96}.

Inverting the Fourier sine transform in \eqref{General solution of PDE and Robin condition} then allows us to write down the solution of the BVP. This is
\begin{align}
\label{Solution Robin case g(t)=0}
\begin{split}
w(S,t)&=\int_0^\infty \varphi (\xi) w^\xi (S,t)d\xi\\
&=\int_0^\infty \frac{2}{\pi}\int_0^{\infty} \zeta(\eta)\sin(\xi \eta) \,
w^\xi (S,t) d \eta d\xi\\
&=\int_0^\infty \int_0^{\infty} \frac{2}{\pi}\zeta(\eta) \sin(\xi \eta) w^\xi (S,t)
d\xi d \eta\\
&=\int_0^\infty \zeta(\eta)\mathcal{G}(S,t,\eta) d \eta,
\end{split}
\end{align}
where on the second line we introduced the inversion integral for the Fourier sine transform.  We have
\begin{align*}
\mathcal{G}(S,\eta,t ) &= \frac{2}{\pi} \int_0^{\infty} \sin(\xi \eta)
w^\xi (S,t)d\xi\\&=\frac{ e^{ct -\frac{\left(\eta - \ln\left( \frac{S}{b} \right)\right)^2}{2 \sigma ^2 t}} S^{\mu}}{\sqrt{2 \pi } \left(\sigma ^2 t\right)^{5/2}} k(S,\eta,t).
\end{align*}
Let $L(S,t)=\gamma   \ln\left( \frac{S}{b} \right) + \tilde{\beta} \sigma ^2 t.$ Then

\begin{align*}
&k(S,\eta,t)=e^{-\frac{2\eta \ln\left( \frac{S}{b} \right)}{ \sigma ^2 t}} \bigg(
 \tilde{\alpha} \sigma ^4 t^2 + \ln\left( \frac{S}{b} \right) \left(2 \gamma  \eta+L(S,t) \right)+\sigma ^2 t \left( \eta \tilde{\beta} -\gamma \right) \\ & + \gamma  \eta ^2 \bigg)
  - \tilde{\alpha} \sigma ^4 t^2 + \ln\left( \frac{S}{b} \right) \left( 2 \gamma  \eta -L(S,t)\right) +\sigma ^2 t \left( \eta \tilde{\beta} + \gamma \right) - \gamma  \eta ^2,\nonumber
\end{align*}
with $\tilde{\alpha}$ and $\tilde{\beta}$ given by \eqref{a tilde} and \eqref{b tilde}.
We would like to write the solution in the form
\begin{align}
    w(S,t)&=\int_0^\infty f(y) p(S,y,t)dy, \nonumber
\end{align}
where $f$ is the initial data and $p(S,y,t)$ is a fundamental solution. In order to do this we use Fubini's Theorem, which leads us to an explicit expression for our fundamental solution  subject to the homogeneous second order Robin boundary conditions. In fact, we can rewrite the solution as
\begin{align}\label{fundamental solution}
        &w(S,t)=\int_0^\infty \zeta(\eta)\mathcal{G}(S,t,\eta)d \eta\nonumber\\
        &=\int_0^\infty
       \frac{2b^{-\mu }}{\tilde{\beta } \sqrt{\Delta}}  \left( \gamma f(b) e^{\frac{\tilde{\beta }\eta}{2 \gamma }} \sinh \left(\frac{\eta \sqrt{\Delta}}{2 \gamma }\right) - \tilde{\beta} \int_0^\eta f\left(be^z \right)\mathcal{K}(\eta,z)
       dz \right)
       \nonumber \\&\quad \times \mathcal{G}(S,t,\eta) d \eta\nonumber\\
	&=\int_0^\infty
       \frac{2b^{-\mu }}{\tilde{\beta } \sqrt{\Delta}} \gamma f(b) e^{\frac{\tilde{\beta }\eta}{2 \gamma }} \sinh \left(\frac{\eta \sqrt{\Delta}}{2 \gamma }\right)
        \mathcal{G}(S,t,\eta) d \eta \nonumber  \\
	& - \int_0^\infty
       \frac{2b^{-\mu }}{\tilde{\beta } \sqrt{\Delta}}  \left( \int_0^\eta f\left(be^z\right)\mathcal{K}(\eta,z)
       dz \right)
        \mathcal{G}(S,t,\eta) d \eta\nonumber\\
& = \ell(S,t)  - \int_0^\infty f\left(be^z \right) \int_z^\infty \frac{2b^{-\mu }}{\tilde{\beta } \sqrt{\Delta}}\mathcal{K}(\eta,z) \mathcal{G}(S,t,\eta) d \eta dz,
\end{align}
where
\begin{equation}
    \ell(S,t):=\int_0^\infty
       \frac{2b^{-\mu }}{\tilde{\beta } \sqrt{\Delta}} \gamma f(b) e^{\frac{\tilde{\beta }\eta}{2 \gamma }} \sinh \left(\frac{\eta \sqrt{\Delta}}{2 \gamma }\right)
        \mathcal{G}(S,t,\eta) d \eta.
\end{equation}
If we set
\begin{equation}
\label{Kernel (transition density)}
	\Bar{p}(S,z,t) := \int_z^{\infty}\frac{2b^{-\mu }}{\tilde{\beta} \sqrt{\Delta}} \mathcal{K}(\eta,z)\mathcal{G}(S,t,\eta)d \eta,
\end{equation}
then we can express  $w(S,t)$ as
\begin{align}
\label{fundamental solution part 2}
    \begin{split}
        w(S,t)& = \ell(S,t)  - \int_0^\infty f\left(be^z\right) \Bar{p} (S,z,t) dz \\
& = \ell(S,t)  - \int_b^\infty f\left(y\right) \tilde{p} (S,y,t) dy, \\
    \end{split}
\end{align}
where
\begin{equation}
\label{kernel of trans dens}
	\tilde{p}(S,y,t) := \frac{1}{y} \Bar{p}\left( S, \ln \left( \frac{y}{b} \right) ,t \right).
\end{equation}
Lastly, formula \eqref{fundamental solution part 2} can be rewritten as
\begin{equation}
    w(S,t) = \int_b^\infty f\left(y\right) p(S,y,t) dy,
\end{equation}
where
\begin{equation}
\label{trans dens up to delta}
	p(S,y,t) = \frac{ \ell(S,t)}{f(b)}\delta^{(b)}(y) - \tilde{p} (S,y,t),
\end{equation}
with $\delta^{(b)}$ denoting the Dirac delta centered in $b \in \R$.
The kernel in \eqref{kernel of trans dens} can be computed explicitly in Mathematica in terms of Gaussians and error functions. It is
complicated so we present it in Appendix 1.

\begin{rem}
We have   solved the second order Robin problem in the case of $g(t)= 0$. We  can also solve  the standard Robin problem by this method.  We sketch the calculation. We begin with \eqnref{General solution of PDE and Robin condition} and apply the boundary condition.  The analysis is very similar to the second order case, though somewhat easier. One finds that a solution satisfying the boundary condition  is
 \begin{align}
        w^\xi (S,t) =e^{-\frac{\xi^2 \sigma^2}{2}t+ct}S^{\mu}\left( A\cos (\xi \ln S )+B \sin (\xi \ln S )\right)
    \end{align}
    where $c=-\frac{1}{2\sigma^2}\left(r-\frac{\sigma^2}{2}\right)^2$ and
\begin{align}
\label{Value of A}
    A=& \beta \xi \cos (\xi \ln b ) +
    b \alpha \sin (\xi \ln b ) + \beta \mu \sin (\xi \ln b )
    \\
    \label{Value of B}
    B=&
    -b \alpha \cos (\xi \ln b ) - \beta \mu \cos (\xi \ln b ) + \beta \xi \sin (\xi \ln b ).
\end{align}
Setting $u(S,t)=\int_0^\infty\varphi(\xi)w^\xi (S,t)d\xi$ and imposing the initial condition leads to an integral equation which also reduces to a Fourier sine transform. The difference is  that in place of the second order equation \eqref{ODE for Boundary Condition} for $\zeta$, we have  $\beta \zeta'(x) + (b \alpha + \beta \mu) \zeta(x) = b^{-\mu} e^{\mu x}f(be^{-x}),$ with $\zeta(0)=0$. This is first order. Obtaining $\zeta$ we proceed exactly as in the second order case and obtain the fundamental solution
\begin{align*}
    \begin{split}
        &p(S,y,t)=-\frac{e^{c t + \frac{\sigma ^2 t
        \tilde{A}^2}{2 \beta ^2}} S^{\mu -\frac{\tilde{A}}{\beta }}}{ \beta  y^{\mu +1+\frac{\tilde{A}}{\beta }} } \tilde{A} b^{\frac{2 \tilde{A}}{\beta }} \text{erfc}\left(\frac{\sigma ^2 t \tilde{A}+\beta \left(\ln \frac{b}{y}+\ln\left(\frac{b}{S}\right) \right)}{\sqrt{2 \beta^2  \sigma^2
       t}}\right)+\\
        & \frac{e^{c t} S^{\mu } y^{-\mu -1}}{\sqrt{2\pi \sigma ^2 t}} \left(b^{\frac{2 \ln \frac{b}{y}+\ln S}{\sigma ^2 t}}+S^{\frac{2 \ln \frac{b}{y} +\ln b}{\sigma ^2 t}}\right) e^{\left(-\frac{\left(\ln \frac{b}{y} +\ln S\right)^2+2 \ln b \ln \frac{b}{y} +\ln ^2 b}{2 \sigma ^2 t}\right)}.
    \end{split}
\end{align*}

\end{rem}
\vskip.2cm
The classical approach to the solution of the problem with $g\neq0$ is found in Friedman's book \cite{Fried64}. For completeness we present it in   Appendix 1. Our purpose now is to develop a new method for solution of these problems that does not require any knowledge of the fundamental solution.

\section{A  novel representation using Hilbert transform methods}\label{DaveHilb}
In this  section we construct a solution to
\begin{align}
\label{GRP}
    w_t &= \frac{1}{2}\sigma^2S^2 w_{SS}+rSw_S, \quad S>b>0 \\
    w(S,0)&=f(S),\label{GRP2}\\
    \alpha w(b,t)+\beta w_S(b,t)&=e^{-\frac12\sigma^2\mu^2t}g(t)\label{GRP3},
\end{align}
which does not require a fundamental solution. We will then extend our method to include second order Robin conditions and a significantly wider
class of PDEs. We assume $\alpha, \beta\neq 0 $ for the remainder of the paper, unless stated otherwise. Choosing one of these constants to be zero reduces the BVP to one of either Dirichlet or Neumann type.

Although explicit solutions can be obtained by our method for many interesting problems, in most cases the representations
that we obtain  will require numerical evaluation. However this is true for every representation and is beyond the scope of our study. We note that there is a considerable literature on the numerical evaluation of Hilbert transforms. We refer the reader to Chapter 14 of \cite{Kin2009a} for an introduction to this topic.

Our approach to the problem uses the {\it even} and {\it odd} Hilbert transforms.
For an exhaustive treatment of the Hilbert transform we refer to King's two volume  treatise, \cite{Kin2009a} and \cite{Kin2009b}.
\subsubsection{The odd and even Hilbert transforms}\label{oeht}
For the reader's convenience we introduce here the material  that we require. The odd and even Hilbert transforms  have some very useful properties which we will exploit. Although these are important operators in their own right, they arose originally as special cases of the classical Hilbert transform.

The various Hilbert transforms are given by  integrals, but are defined only in the principal value sense.
\begin{defn}
We will define the {\it even Hilbert transform} of a function $f:[0,\infty)\to\mathbb{R}$ by the principal value integral
\begin{align}\label{even}
(\mathcal{H}_ef)(x)=\frac{2 x}{\pi}\pint\frac{f(y)}{x^2-y^2}d y,
\end{align}
assuming the principal value integral exists.
Similarly we define the {\it odd Hilbert transform} of a function $f:[0,\infty)\to\mathbb{R}$ by the principal value integral
\begin{align}\label{odd}
(\mathcal{H}_of)(x) =\frac{2 }{\pi}\pint\frac{yf(y)}{x^2-y^2}d y,
\end{align}
again assuming convergence in the principal value sense.
\end{defn}
It is sufficient that $f\in L^2[0,\infty)$ for these transforms to exist. If $f\in L^2[0,\infty)$ then $\mathcal{H}_ef\in L^2[0,\infty)$ and $\mathcal{H}_of \in L^2[0,\infty)$. One can extend the operators to other function spaces, but we will avoid a discussion and simply refer the reader to \cite{Kin2009a}.

\begin{rem}
The Hilbert transform of a suitable function $f$ is defined by
\begin{align}
(\mathcal{H}f)(x)=\frac{1}{\pi}\mathrm{P}\int_{-\infty}^\infty \frac{f(s)}{x-s} d s.
\end{align}
Again if $f\in L^2(\mathbb{R})$ then the Hilbert transform exists and $\mathcal{H}f\in L^2(\mathbb{R})$.

The Hilbert transform has many useful properties, most of which can be found in King's books. An extensive table of transform pairs can be found in \cite{Kin2009b}.
The inverse of the Hilbert transform is simply $-\mathcal{H}.$ That is, $\mathcal{H}^2=-I,$ where $I$ denotes the identity operator.

Now suppose that $f$ is even. Then a simple change of variables in the integral yields
\begin{align}
(\mathcal{H}f)(x)=\frac{2 x}{\pi}\pint\frac{f(y)}{x^2-y^2}d y=(\mathcal{H}_ef)(x).
\end{align}
Conversely if $f$ is odd we obtain
\begin{align}
(\mathcal{H}f)(x)=\frac{2 }{\pi}\pint\frac{yf(y)}{x^2-y^2}d y=(\mathcal{H}_of)(x).
\end{align}
So the odd and even Hilbert transforms can be regarded as special cases of the usual Hilbert transform. This is extremely useful for obtaining
  properties of $\mathcal{H}_e$ and $\mathcal{H}_o.$
\end{rem}

The fundamental relationship between the even and odd transforms is  $$ \mathcal{H}_e \mathcal{H}_o= \mathcal{H}_o \mathcal{H}_e=-I.$$ See \cite{Kin2009a}, page 261.
We also note that if $k(x)=xh(x)$ then
\begin{align}\label{e39}
(\mathcal{H}_e k)(x)=\frac{2 x}{\pi}\pint\frac{yh(y)}{x^2-y^2}d y=x(\mathcal{H}_oh)(x).
\end{align}
This fact will be of importance below. For the second order problem we will require the even Hilbert transform of $y^2(\mathcal{H}_of)(y)$, but this is best presented in context. See equation \ref{ysqHo} below.

The odd and even Hilbert transforms arise because of their connection with the Fourier cosine and sine transforms. If $\mathcal{F}_c$ and $\mathcal{F}_s$ are the cosine and sine transforms respectively, then we have $\mathcal{F}_s^{-1}\mathcal{F}_c=\mathcal{H}_e$ and $\mathcal{F}_c^{-1}\mathcal{F}_s=-\mathcal{H}_o$. A proof can be found in {\cite[p.~259]{Kin2009a}.\footnote{ Note: King uses a slightly different definition of the
sine and cosine transforms, by including a multiplicative factor of $\sqrt{\frac{2}{\pi}}$. Our statements are equivalent to his.}
\subsection{The Solution of the Robin Problem}
We will construct a solution of the problem \eqref{GRP}-\eqref{GRP3}.
We will assume $\alpha+\frac{\mu}{
\beta}\neq 0.$ The case where $\alpha+\frac{\mu}{\beta}=0$ can be handled by a modification of our method.

We will use two linearly independent solutions of the PDE to solve our modified problem. These are
\begin{align}
h_1^\xi(S,t)&=\left(\frac{S}{b}\right)^\mu e^{-\frac12\sigma^2t(\xi^2+\mu^2)}\cos\left(\xi\ln\left(\frac{S}{b}\right)\right),\label{h1}\\
h_2^\xi(S,t)&=\left(\frac{S}{b}\right)^\mu e^{-\frac12\sigma^2t(\xi^2+\mu^2)}\sin\left(\xi\ln\left( \frac{S}{b}\right)\right)\label{h2}.
\end{align}

Notice that our BVP has a slightly different form to that used previously. We pose the problem in this way to avoid a technicality
involving the Laplace transform. This arises in taking the inverse Laplace transform of $e^{as}F(s), a>0$. If $a<0$ the inverse
Laplace transform is simply $f(t+a)H(t+a)$ where $f$ is the inverse Laplace transform of $F$ and $H$ is the Heaviside function.

However if $a$ is positive one has to decide what the inverse transform will be. A natural choice is to insist that $f$ be zero to
the left of the origin, in which case the inverse Laplace transform will be $f(t+a).$ However the issue really requires a discussion of the
Laplace transform as a distribution. So we   solve the modified problem and refer the reader to Laurent Schwartz's treatment
of the Laplace transform within the theory of distributions in the book \cite{LS2008}.


We construct a solution of the PDE of the form
\begin{align}\label{PDEsolrep}
w(S,t)=\int_0^\infty\varphi(\xi)h_1^\xi(S,t)d\xi+\int_0^\infty\psi(\xi)h_2^\xi(S,t)d\xi.
\end{align}
Our aim is to choose $\varphi$ and $\psi$ so that the solutions satisfies the boundary and initial conditions.
Hence we must have
\begin{align}
w(S,0)=\left(\frac{S}{b}\right)^\mu\left[\widehat{\varphi}_c\left(\ln\left(\frac{S}{b}\right)\right)+\widehat{\psi}_s\left(\ln\left(\frac{S}{b}\right)\right)\right]=f(S).
\end{align}
Here $\widehat{\varphi}_c$ and $\widehat{\psi}_s$ are the cosine transform of $\varphi$ and the sine transform of $\psi$ respectively. Setting $y=\ln\left(\frac{S}{b}\right)$ we obtain

\begin{align}
\widehat{\varphi}_c(y)+\widehat{\psi}_s(y)=e^{-\mu y}f\left(be^y\right).
\end{align}
Applying the inverse Fourier cosine transform we obtain the relation
\begin{align}\label{oddH}
\varphi(\xi)&=\frac{2}{\pi}\int_0^\infty e^{-\mu y}f\left(be^y\right)\cos(\xi y)d y-(\mathcal{F}^{-1}_c\mathcal{F}_s\psi)(\xi)\nonumber\\
&=\frac{2}{\pi}\int_0^\infty e^{-\mu y}f\left(be^y\right)\cos(\xi y)d y+(\mathcal{H}_o\psi)(\xi).
\end{align}
Set $F(\xi)=\frac{2}{\pi}\int_0^\infty e^{-\mu y}f\left(be^y\right)\cos(\xi y)d y.$ Then $\varphi(\xi)=F(\xi)+(\mathcal{H}_o\psi)(\xi).$

It is easy to see that
\begin{align}
w(b,t)=\int_0^\infty\varphi(\xi)e^{-\frac12\sigma^2t(\xi^2+\mu^2)}d\xi,
\end{align}
and after some straightforward calculations we obtain
\begin{align}
\alpha w(b,t)+\beta w_S(b,t)&=\left(\alpha+\frac{\mu\beta}{b}\right)\int_0^\infty\varphi(\xi)e^{-\frac12\sigma^2t(\xi^2+\mu^2)}d\xi\nonumber\\&+  \frac{\beta}{b}\int_0^\infty\xi\psi(\xi)e^{-\frac12\sigma^2t(\xi^2+\mu^2)}d\xi = e^{-\frac12\sigma^2\mu^2t}g(t)\nonumber.
\end{align}
Cancelling the factor of $e^{-\frac12\sigma^2\mu^2t}$ and using the substitutions $z=\xi^2$, $s=\frac12\sigma^2t$ gives us
\begin{align*}
\left(\alpha+\frac{\mu\beta}{b}\right)\int_0^\infty\frac{\varphi(\sqrt{z})}{2\sqrt{z}}e^{-zs}d z +  \frac{\beta}{b}\int_0^\infty\frac12\psi(\sqrt{z})e^{-zs}d z= g\left(\frac{2s}{\sigma^2}\right).
\end{align*}
We assume that $\left(\alpha+\frac{\mu\beta}{b}\right)\neq 0.$ The case when $\left(\alpha+\frac{\mu\beta}{b}\right)=0$ reduces to the inversion of a single Laplace transform for $\psi.$  From this and the relation $\varphi(\xi)=F(\xi)+(\mathcal{H}_o\psi)(\xi) $ we immediately obtain $\varphi$ and we have our solution.

Assuming that $\left(\alpha+\frac{\mu\beta}{b}\right)\neq 0 $, taking the inverse Laplace Transform and letting $z=\xi^2$ we obtain
\begin{align}\label{e50}
\left(\alpha+\frac{\mu\beta}{b}\right)\varphi(\xi)=2\xi\mathcal{L}^{-1}\left[g\left(\frac{2s}{\sigma^2}\right);\xi^2 \right]-\frac{\xi\beta}{b}\psi(\xi).
\end{align}
Set $K(\xi)=2\xi\mathcal{L}^{-1}\left[ g\left(\frac{2s}{\sigma^2}\right);\xi^2 \right].$ Obviously in the homogeneous case $K(\xi)=0$, so there is no Laplace transform inversion necessary. Using \eqnref{oddH} we have the relation
\begin{align}\label{e42}
\left(\alpha+\frac{\mu\beta}{b}\right)(F(\xi)+ (\mathcal{H}_o\psi)(\xi))=K(\xi)-\frac{\xi\beta}{b}\psi(\xi).
\end{align}

We now take the even Hilbert transform of both sides to obtain
\begin{align}\label{e51}
\left(\alpha+\frac{\mu\beta}{b}\right)\left[(\mathcal{H}_e F)(\xi)-\psi(\xi)\right]=(\mathcal{H}_eK)(\xi)-\frac{\beta}{b}\xi(\mathcal{H}_o\psi)(\xi),
\end{align}
where we used the relation \eqnref{e39}. So we have obtained a pair of simultaneous equations for $\mathcal{H}_o\psi$ and $\psi$.  Clearly equation \eqref{e42} gives us
\begin{align}
(\mathcal{H}_o\psi)(\xi)=\left(\alpha+\frac{\mu\beta}{b}\right)^{-1}\left[K(\xi)-\frac{\xi\beta}{b}\psi(\xi)\right]-F(\xi).
\end{align}
So that
\begin{align}
 (\mathcal{H}_e F)(\xi)-\psi(\xi) =\frac{(\mathcal{H}_eK)(\xi)-\frac{\beta}{b}\xi\left( \frac{\left[K(\xi)-\frac{\xi\beta}{b}\psi(\xi)\right]}{\left(\alpha+\frac{\mu\beta}{b}\right) }-F(\xi)\right)}{\left(\alpha+\frac{\mu\beta}{b}\right)}.
\end{align}

Rearranging this gives

\begin{align*}
\psi(\xi)=\frac{\nu\left[\frac{b\alpha+\mu\beta}{b}(\mathcal{H}_eF)(\xi)-(\mathcal{H}_eK)(\xi)+
\frac{\beta\xi}{b}\left(\frac{b}{\alpha b+\mu\beta} K(\xi)-F(\xi)\right)\right]}{(\alpha b+\mu\beta)^2+\beta^2\xi^2},
\end{align*}
where $\nu=b(\alpha b+\mu\beta)$.
From this we find $\varphi$ and we have obtained a potential solution to the boundary value problem. This leads us to the following result.
\begin{thm}\label{HilbThm}
Let $f\in L^1\left((b,\infty),\frac{d \eta}{\eta^{\mu+1}}\right)$ and $g(t)=\int_0^\infty G(x)e^{-xt}d x$ where $G$, $G'(x)$ and $G''$ are integrable.
Let  $h_1^\xi$ and $h_2^\xi$ be given by \eqref{h1} and \eqref{h2} respectively.
Then if $\alpha+\frac{\mu\beta}{b}\neq 0$ and $\nu=b(\alpha b+\mu\beta)$, the problem \eqref{GRP}-\eqref{GRP3} has a solution given by
\begin{align}\label{soldeed}
w(S,t)=\int_0^\infty\varphi(\xi)h_1^\xi(S,t)d\xi+\int_0^\infty\psi(\xi)h_2^\xi(S,t)d\xi,
\end{align}
 where
\begin{align*}
\psi(\xi)=\frac{\nu\left[\frac{b\alpha+\mu\beta}{b}(\mathcal{H}_eF)(\xi)-(\mathcal{H}_eK)(\xi)+
\frac{\beta\xi}{b}\left(\frac{b}{\alpha b+\mu\beta} K(\xi)-F(\xi)\right)\right]}{(\alpha b+\mu\beta)^2+\beta^2\xi^2},
\end{align*}
and
\begin{align*}
\varphi(\xi)=\left(\alpha+\frac{\mu\beta}{b}\right)^{-1}\left[\sigma^2\xi G\left(\frac{\sigma^2\xi^2}{2}\right)-\frac{\beta \xi}{b}\psi(\xi)\right].
\end{align*}
Here $F(\xi)=\frac{2}{\pi}\int_0^\infty e^{-\mu y}f\left(be^y\right)\cos(\xi y)d y $ and $K(\xi)=\sigma^2\xi G\left(\frac{\sigma^2\xi^2}{2}\right).$
\end{thm}
\begin{proof}
To complete the proof, we must establish sufficient conditions to guarantee that we do in fact have a solution. To this end, we make the following observations. Since $\mathcal{H}_e=\mathcal{F}_s^{-1}\mathcal{F}_c$ we easily obtain
\begin{align*}
(\mathcal{H}_eF)(\xi)&=\mathcal{F}_s^{-1}\mathcal{F}_c\frac{2}{\pi}\int_0^\infty e^{-\mu y}f\left(be^y\right)\cos(\xi y)d y
\\&=\frac{2}{\pi}\int_0^\infty e^{-\mu y}f(be^y)\sin(\xi y)d y\\
&=\frac{2}{\pi}\int_b^\infty \left(\frac{\eta}{b}\right)^{-\mu}f(\eta)\sin\left(\xi\ln\left(\frac{\eta}{b}\right)\right)\frac{d \eta}{\eta}.
\end{align*}

Let $f\in L^1((b,\infty),\eta^{-\mu-1}d \eta)$ and denote $\|k\|_b=\int_b^\infty|k(x)|x^{\mu-1}d x.$ Then we obtain the inequality
\begin{align}
|\mathcal{H}_eF|\leq \frac{2}{\pi b^{\mu}}\|f\|_b.
\end{align}
It is also easy to see that $|F|\leq \frac{2}{\pi b^{\mu}}\|f\|_b.$ From this we can conclude that
\begin{align}\label{deriv1}
\frac{\partial^j}{\partial S^j}\left(\frac{-b(\alpha b+\mu\beta) }{(\alpha b+\mu\beta)^2+\beta^2\xi^2}(\mathcal{H}_eF)(\xi)h_k^{\xi}(S,t)\right),
\end{align}
is integrable for each $k=1,2$, $j=0,1,2$ . Similarly
\begin{align}\label{deriv2}
\frac{\partial }{\partial t}\left(\frac{-b(\alpha b+\mu\beta) }{(\alpha b+\mu\beta)^2+\beta^2\xi^2}(\mathcal{H}_eF)(\xi)h_k^{\xi}(S,t)\right),
\end{align}
is integrable. Both facts follow from the Gaussian decay of the solutions $h_k^{\xi}(S,t).$ Specifically we can bound \eqref{deriv1} in $\xi$ by
\begin{align}
\left|\frac{\partial^j}{\partial S^j}\left(\frac{-b(\alpha b+\mu\beta) }{(\alpha b+\mu\beta)^2+\beta^2\xi^2}(\mathcal{H}_eF)(\xi)h_k^{\xi}(S,t)\right)\right|\leq C \xi^j e^{-\gamma \xi^2},
\end{align}
for some positive constants $\gamma$, $C$, which will depend on $S$. Similarly for (\ref{deriv2}).

The same argument shows that
\begin{align*}
\frac{\partial^j}{\partial S^j}\left(\frac{  b\beta\xi}{(\alpha b+\mu\beta)^2+\beta^2\xi^2}F(\xi)h_k^{\xi}(S,t)\right),
\end{align*}
and
\begin{align*}
\frac{\partial }{\partial t}\left(\frac{-b\beta\xi}{(\alpha b+\mu\beta)^2+\beta^2\xi^2}F(\xi)h_k^{\xi}(S,t)\right),
\end{align*}
are also integrable for $k=1,2$ and $j=0,1,2$.

Next we suppose that $g(t)=\int_0^\infty G(x)e^{-xt}d x$ and for simplicity we will suppose that $\|G\|_{L^1}=\int_0^\infty |G(x)|d x<\infty.$
In fact this assumption can be relaxed to allow $\int_0^\infty|e^{rx}G(x)|d x<\infty$ for some $r>0$. We can also allow $G$ to be a distribution.
We will not go into these technicalities here, but we will present an example when $G$ is a distribution below.

It immediately follows from our assumption that $$K(\xi)=2\xi\mathcal{L}^{-1}\left[ g\left(\frac{2s}{\sigma^2}\right);\xi^2 \right]=\sigma^2\xi G\left(\frac{\sigma^2\xi^2}{2}\right).$$
Consequently a simple change of variables gives $\|K\|_{L^1}=\|G\|_{L^1}<\infty.$
Arguing as previously we see that

\begin{align*}
\frac{\partial^j}{\partial S^j}\left(\frac{b\beta \xi}{(\alpha b+\mu\beta)^2+\beta^2\xi^2}K(\xi) h_k^{\xi}(S,t)\right),
\end{align*}
for $k=1,2$, $j=0,1,2$ are integrable, as is
\begin{align*}
\frac{\partial }{\partial t}\left(\frac{b\beta \xi}{(\alpha b+\mu\beta)^2+\beta^2\xi^2}K(\xi) h_k^{\xi}(S,t)\right).
\end{align*}

Turning to $\mathcal{H}_eK $ we observe that if $k$ is twice differentiable and $k',k''$ are integrable, then the Fourier cosine transform $\mathcal{F}_ck$ is also integrable. This is a famous result and follows via integration by parts, see \cite{YK68}.

The observation that $|\mathcal{F}_c k|\leq \|k\|_{L^1} $ is elementary. So we have
$$ |\mathcal{F}_s^{-1}(\mathcal{F}_c k)|\leq \frac{2}{\pi} |\mathcal{F}_ck | \leq \frac{2}{\pi}\|k\|_{L^1}. $$
Hence we can conclude that if $G'$ and $G''$ are integrable, then $|\mathcal{H}_e K|\leq \frac{2}{\pi}\|G\|_{L^1}.$

This implies that for each $k=1,2$, $j=0,1,2$

\begin{align*}
\frac{\partial^j}{\partial S^j}\left( \frac{ (\alpha b+\mu\beta)^2}{(\alpha b+\mu\beta)^2+\beta^2\xi^2}(\mathcal{H}_eK)(\xi)h_k^{\xi}(S,t)\right),
\end{align*}
is integrable, as is
\begin{align*}
\frac{\partial }{\partial t}\left( \frac{ (\alpha b+\mu\beta)^2}{(\alpha b+\mu\beta)^2+\beta^2\xi^2}(\mathcal{H}_eK)(\xi)h_k^{\xi}(S,t)\right).
\end{align*}

From this we see that if
\begin{align}
w(S,t)=\int_0^\infty\varphi(\xi)h_1^\xi(S,t)d\xi+\int_0^\infty\psi(\xi)h_2^\xi(S,t)d\xi,
\end{align}
then three applications of the Dominated Convergence Theorem allows us to differentiate under the integral sign with respect to $t$ and twice with respect
to $S$. Since $h_1^\xi(S,t)$ and $h_2^\xi(S,t)$ satisfy the PDE, it follows that
\eqnref{soldeed} also solves the PDE. By construction it also satisfies the boundary conditions, as well as the initial data.

\end{proof}

\begin{rem} Thus we have constructed a solution of the BVP that does not require us to known a fundamental solution for
the problem. An alternative representation for the solution of a BVP can be very useful. It may be more tractable or more computationally efficient than the
classical method. However a discussion of this issue is beyond the scope of the current work. We content ourselves with an illustrative example.
\end{rem}

\begin{example}
We take $g(t)=t $, $f(S)= \frac{1}{\sigma^2}\left(\alpha+\frac{\mu\beta}{b}\right)^{-1}\left(\frac{S}{b}\right)^\mu\left(\ln\left(\frac{S}{b}\right)\right)^2 $ and assume that $\alpha+\frac{\mu\beta}{b}\neq0.$ Now recall that if $\delta$ is the Dirac delta function then by definition,
\begin{align}
\int_0^\infty\delta'(x)f(x)d x=-f'(0).
\end{align}
Hence $\int_0^\infty \delta'(x)e^{-sx}d x=s.$ After some calculations that we omit, it turns out that we can take $\psi=0$ and
\begin{align}
\varphi(\xi)= \left(\alpha+\frac{\mu\beta}{b}\right)^{-1}\frac{4\xi}{\sigma^2}\delta'(\xi^2).
\end{align}
Consequently a solution of our problem with $g(t)=t$ is
\begin{align*}
h(S,t)&=\int_0^\infty  \frac{4\xi}{\sigma^2\left(\alpha+\frac{\mu\beta}{b}\right)}\delta'(\xi^2)\left(\frac{S}{b}\right)^\mu e^{-\frac12\sigma^2t(\xi^2+\mu^2)}\cos\left(\xi\ln \frac{S}{b}\right)d\xi\nonumber\\
&=\frac{2\left(\frac{S}{b}\right)^\mu }{\sigma^2\left(\alpha+\frac{\mu\beta}{b}\right)} e^{-\frac12\sigma^2\mu^2 t } \int_0^\infty \delta'(z) e^{-\frac12\sigma^2z t }\cos\left(\sqrt{z}\ln \frac{S}{b}\right)d z\nonumber \\
&=\frac{2\left(\frac{S}{b}\right)^\mu}{\sigma^2\left(\alpha+\frac{\mu\beta}{b}\right)} e^{-\frac12\sigma^2\mu^2 t }\left[\frac{\sigma^2 }{2}+\frac12\left(\ln\left(\frac{S}{b}\right)\right)^2\right].
\end{align*}
It is not hard to see that
\begin{align*}
h(S,0)
&=\frac{2}{\sigma^2}\left(\alpha+\frac{\mu\beta}{b}\right)^{-1}\left(\frac{S}{b}\right)^\mu \int_0^\infty \delta'(z) \cos\left(\sqrt{z}\ln \frac{S}{b}\right)d z\\
&=\frac{1}{\sigma^2}\left(\alpha+\frac{\mu\beta}{b}\right)^{-1}\left(\frac{S}{b}\right)^\mu\left(\ln\left(\frac{S}{b}\right)\right)^2 .
\end{align*}
So the initial condition is satisfied. Using Mathematica it is easy to see that the boundary conditions are satisfied. Choosing $g$ to be a polynomial
leads to the appearance of derivatives of the Dirac delta and integrals involving these distributions are particularly easy to evaluate.
\end{example}

\section{Hilbert Transform Methods for Second Order Boundary Conditions}\label{DaveHilb2}
The result of the previous section give rise to a number of questions. The first is whether the Hilbert transform approach can be extended to the second
order Robin problem? A more important question is whether the method can be applied to boundary value problems for other PDEs? The answer to both
questions is yes.

We start with the extension to the second order Robin problem, then turn to the second question. We will solve
\begin{align}
\label{P3}
    w_t &= \frac{1}{2}\sigma^2S^2 w_{SS}+rSw_S, \\
    w(S,0)&=f(S),\label{P4}\\
    \alpha w(b,t)+\beta w_S(b,t)+\gamma w_{SS}(b,t)&=e^{-\frac12\sigma^2\mu^2t}g(t),\label{P5}
\end{align} where $  S>b>0. $

In fact the same method can be used to solve this problem, but with an extra step. Our solution will again be of the form (\ref{PDEsolrep}) as in the previous section. It is straightforward to show that $w(b,0)=\int_0^\infty\varphi(\xi)d\xi=f(b).$ This fact will be useful.

As in the regular Robin condition case we deduce that $\varphi(\xi)=F(\xi)+(\mathcal{H}_o\psi)(\xi),$
where $F$ is exactly as before. Here we make the observation that $$\int_0^\infty (\mathcal{H}_o\psi)(\xi)d\xi=\int_0^\infty(\varphi(\xi)-F(\xi))d\xi=f(b)-\int_0^\infty F(\xi)d\xi.$$
We here assume that $F$ is integrable. This can be guaranteed by imposing mild conditions on $f$. We will give a sufficient condition for integrability  below.

Introduce the constants $A=\alpha+\frac{1}{b^2}\gamma\mu(\mu-1)+\frac{\mu\beta}{b}$ and $B=\frac{\beta}{b}+\frac{\gamma}{b^2}(2\mu-1).$ Then the boundary condition yields the equation
\begin{align}
 \int_0^\infty(A-\frac{\gamma}{b^2}\xi^2) \varphi(\xi)e^{-\frac12\sigma^2t \xi^2 }d\xi +B\int_0^\infty
\xi\psi(\xi)e^{-\frac12\sigma^2t \xi^2 }d\xi= g(t).
\end{align}
Converting the integrals to Laplace transforms and inverting as previously we obtain
\begin{align*}
(A-\frac{\gamma}{b^2}\xi^2)\varphi(\xi)+B\xi\psi(\xi)=2\xi\mathcal{L}^{-1}\left[g\left(\frac{2s}{\sigma^2}\right);\xi^2 \right].
\end{align*}
We will insist that $A$ and $\gamma$ have opposite signs to ensure that $\varphi$ is nonsingular. Let $G(\xi)=2\xi\mathcal{L}^{-1}\left[g\left(\frac{2s}{\sigma^2}\right);\xi^2 \right].$ Then we can rewrite the equation as
\begin{align*}
(A-\frac{\gamma}{b^2}\xi^2)[F(\xi)+(\mathcal{H}_0\psi)(\xi)]+B\xi\psi(\xi)=G(\xi).
\end{align*}
This is the same as
\begin{align}\label{psieqn2}
A(\mathcal{H}_0\psi)(\xi)-\frac{\gamma}{b^2}\xi^2(\mathcal{H}_0\psi)(\xi)+B\xi\psi(\xi)=G(\xi)-(A-\frac{\gamma}{b^2}\xi^2) F(\xi).
\end{align}

To proceed we need to take the even Hilbert transform of both sides.  We make the observation that
\begin{align}
\mathcal{H}_e(y^2(\mathcal{H}_0\psi)(y))(x)&=\frac{2 x}{\pi}\pint\frac{y^2(\mathcal{H}_0\psi)(y)}{x^2-y^2}d y\nonumber\\
&=\frac{2 x}{\pi}\pint\frac{(y^2-x^2+x^2)(\mathcal{H}_0\psi)(y)}{x^2-y^2}d y\nonumber\\
&=-\frac{2x}{\pi}\pint (\mathcal{H}_0\psi)(y)d y+x^2(\mathcal{H}_e\mathcal{H}_o\psi)(x)\nonumber\\
&=-\frac{2x}{\pi}\left(f(b)-\pint F(y)d y\right)-x^2\psi(x)\nonumber\\
&=Cx-x^2\psi(x),\label{ysqHo}
\end{align}
where we have defined the constant $C=-\frac{2 }{\pi}\left(f(b)-\int_0^\infty F(y)d y\right).$ We dropped the principal value because we are assuming
suitable integrability. We can of course weaken this assumption and reintroduce the principal values of the integrals.

Taking the even Hilbert transform of (\ref{psieqn2}) gives us the relation
\begin{align}\label{oma}
-(A-\frac{\gamma}{b^2}\xi^2) \psi(\xi)- \frac{\gamma}{b^2} C\xi +B\xi(\mathcal{H}_0\psi)(\xi)=\mathcal{H}_e[G(\xi)-(A-\frac{\gamma}{b^2}\xi^2) F(\xi)].
\end{align}
Now we easily see that
\begin{align*}
(\mathcal{H}_e x^2F(x))(\xi)=-\frac{2\xi}{\pi}\pint F(x)d x+\xi^2(\mathcal{H}_eF)(\xi).
\end{align*}
So our assumptions on $F$ and the assumptions of the previous section are sufficient to guarantee that the right hand side of (\ref{oma}) exists.

We know that $\varphi(\xi)=\frac{G(\xi)- B\xi\psi(\xi)}{A-\frac{\gamma}{b^2}\xi^2},$
so that
\begin{align*}
(\mathcal{H}_0\psi)(\xi)=\frac{G(\xi)- B\xi\psi(\xi)}{A-\frac{\gamma}{b^2}\xi^2}-F(\xi).
\end{align*}

We thus have
\begin{align*}
-A\psi(\xi)- \frac{\gamma}{b^2}(C\xi-\xi^2\psi(\xi))+B\xi\left(\frac{G(\xi)- B\xi\psi(\xi)}{A-\frac{\gamma}{b^2}\xi^2}-F(\xi)\right)=M(\xi),
\end{align*}
with $M(\xi)=\mathcal{H}_e[G(\xi)-(A-\frac{\gamma}{b^2}\xi^2) F(\xi)].$

This gives
\begin{align*}
-\left(A-\frac{\gamma}{b^2}\xi^2+\frac{B^2\xi^2}{A-\frac{\gamma}{b^2}\xi^2}\right)\psi=C\frac{\gamma}{b^2}\xi+M(\xi)+B\xi F(\xi)-B\xi \frac{G(\xi) }{A-\frac{\gamma}{b^2}\xi^2}.
\end{align*}

We let
\begin{align}\label{N}
N(\xi)=C\frac{\gamma}{b^2}\xi+M(\xi)+B\xi F(\xi)-B\xi \frac{G(\xi) }{A-\frac{\gamma}{b^2}\xi^2}.
\end{align} Then we have
\begin{align}\label{2psiform}
\psi(\xi)=\frac{-N(\xi)(A-\frac{\gamma}{b^2}\xi^2)}{ (A-\frac{\gamma}{b^2}\xi^2)^2+ B^2\xi^2  }.
\end{align}
Combining this we have the following result.
\begin{thm}\label{HilbThm2} We suppose that $f$ and $g$ satisfy the same conditions as in Theorem \ref{HilbThm}. Suppose further that $\eta f'(\eta),\eta^2f''(\eta)\in L^1\left((b,\infty),\frac{d \eta}{\eta^{\mu+1}}\right)$. Let $F,G$ be as in Theorem \ref{HilbThm}. Then the Problem \ref{P3}, with $\alpha,\beta,\gamma\neq 0$, has a solution given by
\begin{align}\label{PDEsolrep2}
w(S,t)=\int_0^\infty\varphi(\xi)h_1^\xi(S,t)d\xi+\int_0^\infty\psi(\xi)h_2^\xi(S,t)d\xi,
\end{align}
where $h_{1,2}^\xi(S,t)$ are as in Theorem \ref{HilbThm} and
\begin{align}\
\psi(\xi)=\frac{-N(\xi)(A-\frac{\gamma}{b^2}\xi^2)}{ (A-\frac{\gamma}{b^2}\xi^2)^2+ B^2\xi^2  }
\end{align}
and \begin{align*}
\varphi(\xi)=\frac{G(\xi)- B\xi\psi(\xi)}{A-\frac{\gamma}{b^2}\xi^2}.
\end{align*}
Here $ A=\alpha+\frac{1}{b^2}\gamma\mu(\mu-1)+\frac{\mu\beta}{b}$, $B=\frac{\beta}{b}+\frac{\gamma}{b^2}(2\mu-1),$  $N$ is as defined in \eqref{N}  and $A$ and $\gamma$ have opposite signs.
\end{thm}
\begin{proof}
The proof follows along the lines of the proof of Theorem \ref{DaveHilb}. So here we only consider the integrability of $F$. We recall that $$F(\xi)=\frac{2}{\pi}\int_0^\infty e^{-\mu y}f\left(be^y\right)\cos(\xi y)d y.$$ Obviously we can write
\begin{align*}
\int_0^\infty F(\xi)d\xi=\int_0^1 F(\xi)d\xi+\int_1^\infty F(\xi)d\xi.
\end{align*}
A simple application of Fubini's Theorem shows that $|\int_0^1 F(\xi)d\xi|<\infty$ provided that $$\int_0^1\left|e^{-\mu y}f\left(be^y\right)\right|\frac{d y}{y}<\infty.$$
The obvious change of variables gives us
\begin{align*}
\int_0^1\left|e^{-\mu y}f\left(be^y\right)\right|\frac{d y}{y}=b^\mu\int_b^{be}|f(\eta)|\frac{d \eta}{\eta^{\mu+1}}.
\end{align*}
Thus $f\in L^1((b,\infty),\frac{d \eta}{\eta^{\mu+1}} )$  guarantees that $|\int_0^1 F(\xi)d\xi|<\infty$.

Now we assume that $f$ is twice differentiable. Integration by parts gives us
\begin{align*}
\int_1^\infty\int_0^\infty e^{-\mu y}f\left(be^y\right) \cos(\xi y)d y &=\int_1^\infty\int_0^\infty \frac{d}{d y}\left(e^{-\mu y}f\left(be^y\right)\right) \frac{\sin(\xi y)}{-\xi}d y d\xi\\
&=\int_1^\infty\left(\left[\frac{d}{d y}(e^{-\mu y}f\left(be^y\right))\frac{\cos(\xi y)}{\xi^2}\right]_0^\infty\right.\\
& \left.-\int_0^\infty \frac{d^2}{dy^2}(e^{-\mu y}f\left(be^y\right) )\frac{\cos(\xi y)}{\xi^2}d y\right)d\xi\\
&=\int_1^\infty \frac{\mu f(b)-bf'(b)}{\xi^2}d\xi -I.
\end{align*}
where $I=\int_1^\infty \int_0^\infty \frac{d^2}{d y^2}(e^{-\mu y}f\left(be^y\right) )\frac{\cos(\xi y)}{\xi^2}d y d \xi.$ Clearly $$\int_1^\infty \frac{\mu f(b)-bf'(b)}{\xi^2}d\xi=\mu f(b)-bf'(b).$$ Further
\begin{align*}
|I|&=\left|\int_1^\infty \int_0^\infty \frac{d^2}{dy^2}(e^{-\mu y}f\left(be^y\right) )\frac{\cos(\xi y)}{\xi^2}d y d\xi\right|\\&\leq \int_1^\infty \int_0^\infty \left|\frac{d^2}{d y^2}(e^{-\mu y}f\left(be^y\right) )\frac{\cos(\xi y)}{\xi^2}\right|d y d\xi\\
&\leq \int_0^\infty \left|\frac{d^2}{d y^2}(e^{-\mu y}f\left(be^y\right) ) \right|d y <\infty,
\end{align*}
where we have employed Fubini's Theorem on the assumption that the final integral is finite. Thus the integrability of $\frac{d^2}{d y^2}(e^{-\mu y}f\left(be^y\right) )$ is sufficient to guarantee that $|\int_0^\infty F(\xi)d\xi|<\infty.$
Since
\begin{align*}
\int_0^\infty \left|\frac{d^2}{d y^2}(e^{-\mu y}f\left(be^y\right) ) \right|d y&=b^\mu\int_b^\infty \left|(Df)(\xi)\right|\frac{d \eta}{\eta^{\mu+1}},
\end{align*}
$(Df)(\xi)=\mu^2f(\eta)+b(1-2\mu)\eta f'(\eta)+\mu^2\eta^2f''(\eta)$,
then it is sufficient to require $f(\eta),\eta f'(\eta),\eta^2f''(\eta)\in L^1((b,\infty),\frac{d \eta}{\eta^{\mu+1}} )$, in order for the integral of $F$ to exist. This also justifies our use of Fubini's Theorem.

The proof of the integrability of $\varphi(\xi)h_1^\xi(S,t)$ and $\psi(\xi)h_2^\xi(S,t)$ as well as the differentiability of the integrals defining the solution proceeds along the same lines as in the proof of Theorem \ref{HilbThm}.

\end{proof}

Again we have a representation of the solution of the second order Robin problem, that does not rely upon any knowledge of a fundamental solution. Given the complexity of the fundamental solution in the second order case, this may be preferred. The problem of efficient evaluation of the solutions will be treated elsewhere. However we remark  that the even and odd Hilbert transforms can be written in terms of the Fourier sine and cosine transforms. Hence
techniques for the evaluation of these transforms can be employed to calculate the even and odd Hilbert transforms. In particular one has access to   Fast Fourier Transform  methods.

\section{A General Result For Second Order Robin Problems}\label{QE2d}

It is possible to construct general formulae for different types of boundary conditions. We present one particular case for
second order Robin conditions. In the simplest version of the theory we need a family of solutions of the form
\begin{align}
w_1(x,t;\xi)=\rho(x)\cos(\xi\mu(x))e^{-\xi^2t},\label{f1}\\
w_2(x,t;\xi)=\rho(x)\sin(\xi\mu(x))e^{-\xi^2t}\label{f2},
\end{align}
where $\mu$ and $\rho$ are assumed to be twice continuously differentiable. There are many variations on this and we will present some examples below. By linearity, if $w_1$ and $w_2$ are solutions of a linear PDE, then so are
\begin{align}
\bar{w}_1(x,t;\xi)=\rho(x)\cos(\xi(\mu(x)-\mu(a)))e^{-\xi^2t},\label{f3}\\
\bar{w}_2(x,t;\xi)=\rho(x)\sin(\xi(\mu(x)-\mu(a)))e^{-\xi^2t}\label{f4}.
\end{align}
In fact many PDEs on the line have solutions of this form. A simple characterisation can be obtained as follows.

If $w_1$ is a solution of the PDE $u_t=P(x)u_{xx}+Q(x)u_x+R(x)u$, where $P,Q,R$ are smooth functions, then we must have
\begin{align*}
&P(x)(\rho''(x)\cos(\xi\mu(x))-2\xi\rho'(x)\sin(\xi\mu(x))-\xi^2\rho(x)\cos(\xi\mu(x)))\\&+   Q(x)(\rho'(x)\cos(\xi\mu(x))-\xi \rho(x)\mu'(x)\sin(\xi\mu(x))\\&+R(x)(\rho(x)\cos(\xi\mu(x))+\xi^2)=0.
\end{align*}

Substitution into the PDE shows that in order for this to be a solution we require $P,Q,R,\mu$ and $\rho$ to satisfy the following system of
equations.
\begin{align}
P(x)\rho''(x)+Q(x)\rho(x)+R(x)\rho(x)&=0,\label{e1} \\
1-P(x)(\mu'(x))^2&=0,\label{e2}\\
Q(x)\rho(x)\mu'(x)+2P(x)\mu'(x)\rho'(x)+P(x)\rho(x)\mu''(x)&=0\label{e3}.
\end{align}
This implies that $\mu'(x)=\pm(P(x))^{-1/2}.$

Plainly $\rho$ must be a time independent solution of the PDE.
It is not hard to show that $w_2(x,t;\xi)$ will be a solution of the
PDE whenever $w_1$ is and vice versa. Most PDEs will not possess solutions of this form, but there are numerous PDEs which do. These defining equations allow us to generate  examples and we will present some below.

Now we present a general result for the solution of second order Robin problems, where the PDE has solutions of this form. We observe that by taking $\gamma\to 0$ we obtain the solution of the standard Robin problem. There are various special cases that we will not cover. They can all be handled by obvious modifications of our arguments.  We make several simplifying assumptions to make the proof easier. More precise conditions are possible, but we will not present a full analysis here. The actual calculations are essentially the same as in the Black-Scholes case, so we leave these to the reader.

\begin{thm}\label{AGentheorem}
Suppose that the PDE
\begin{align}\label{genrobthm}
u_t=P(x)u_{xx}+Q(x)u_x+R(x)u, \ x>b, t>0,
\end{align}
has solutions given by \eqref{f3} and \eqref{f4}. Suppose that $\mu:[b,\infty)\to\mathbb{R}$ is invertible, $g$ has an inverse Laplace transform and
$\tilde{f}\in L^1[0,\infty)$ where $\tilde{f}(y)=\frac{f(\mu^{-1}(y+\mu(b)))}{\rho(\mu^{-1}(y+\mu(b)))}.$

Define the constants $A =\alpha\rho(b)+\beta\rho'(b)+\gamma\rho''(b),
B =-\gamma\rho(b)(\mu'(b))^2,$
$C=\beta\rho(b)\mu'(b)+\gamma(\rho(b)\mu''(b)+2\rho'(b)\mu'(b)),$
$\chi =\frac{2B}{\pi}\left(c-\int_0^\infty F(\xi)d\xi\right),$
and  $c=\lim_{z\to b}\frac{f(z)}{\rho(z)}$, which we assume exists. Suppose that $A$ and $B$ are non zero and have opposite signs.

Let $F(\xi)=\frac{2}{\pi}\int_0^\infty \tilde{f}(\eta)\cos(\xi\eta)d\eta$,
$G(\xi)=2\xi\mathcal{L}^{-1}[g(t);\xi^2]$ and suppose that $K(\xi)=G(\xi)-(A-B\xi^2)F(\xi)\in L^2[0,\infty).$
Then a solution of \eqref{genrobthm} satisfying $u(x,0)=f(x)$ and $\alpha u(b,t)+\beta u_x(b,t)+\gamma u_{xx}(b,t)=g(t)$
is
\begin{align}\label{solution}
u(x,t)&=\int_0^\infty \varphi(\xi)\rho(x)\cos(\xi(\mu(x)-\mu(b)))e^{-\xi^2t}d\xi\nonumber \\&+\int_0^\infty \varphi(\xi)\rho(x)\cos(\xi(\mu(x)-\mu(b)))e^{-\xi^2t}d\xi,
\end{align}
in which
\begin{align}
\psi(\xi)&=-\frac{(A-B\xi^2)((\mathcal{H}_e K)(\xi)-\chi\xi)-C\xi K(\xi)}{(A-B\xi^2)^2+C^2\xi^2},
\end{align}
and
\begin{align}
\varphi(\xi)=\frac{G(\xi)-C\xi\psi(\xi)}{A-B\xi^2}.
\end{align}

\end{thm}

\begin{proof}
The derivation of the solutions follows the same line as the Black-Scholes case. The initial condition gives
$\varphi(\xi)=F(\xi)+(\mathcal{H}_o\psi)(\xi)$
and the boundary condition gives us $[A-B\xi^2]\varphi(\xi)+C\xi\psi(\xi)=G(\xi).$ We immediately obtain
\begin{align*}
[A-B\xi^2]\psi(\xi)+C\xi\psi(\xi)=K(\xi).
\end{align*}
Taking the even Hilbert transform of both sides gives us
\begin{align}
-(A-B\xi^2)\psi(\xi)+C\xi(\mathcal{H}_o)\psi)(\xi)+\chi\xi=(\mathcal{H}_e K)(\xi).
\end{align}
We therefore have a pair of simultaneous equations. Solving gives us our
expression for $\psi$ and $\varphi.$

Now $w_1(x,t;\xi)$ and $w_2(x,t;\xi)$ are solutions of the PDE and they have Gaussian decay in $\xi.$ By Riesz's inequality
(see \cite{Kin2009a}), there exists a constant $\mathcal{R}$ such that $\|\mathcal{H}_eK\|_2\leq \mathcal{R}\|K\|_2.$
It immediately follows from H\"{o}lder's inequality and our assumptions on $F$ and $G$ that the integral defining the solution is convergent.
Similarly the integrals
\begin{align*}
\int_0^\infty \eta(\xi)\frac{\partial^j}{\partial x^j}w_{1,2}(x,t;\xi)d\xi,
\end{align*}
$j=1,2,$ and
\begin{align*}
\int_0^\infty \eta(\xi)\frac{\partial }{\partial t}w_{1,2}(x,t;\xi)d\xi,
\end{align*}
are convergent. Here $\eta$ is either $\varphi$ or $\psi.$ Therefore \eqref{solution} is a solution of the PDE and by construction it
satisfies the boundary conditions and the initial condition.
\end{proof}

\subsection{Some Examples with Elementary Solutions}\label{exes}
Using the previous remarks, one can generate an interesting variety of examples.
\begin{example} For the heat equation with drift
\begin{align}
u_t=u_{xx}+au_x,\ x>0, t>0,
\end{align}
we let  $u(x,0)=f(x)$, $\alpha u(0,t)+\beta u_x(0,t)+\gamma u_{xx}(0,t)=e^{-\frac{a^2}{4}t}g(t).$ The elementary solutions are $w_1(x,t;\xi)=e^{-\xi^2t-\frac{a^2}{4}t}e^{-\frac{a}{2}x}\cos(\xi x),  w_2(x,t;\xi)=e^{-\xi^2t-\frac{a^2}{4}t}e^{-\frac{a}{2}x}\sin(\xi x)$. Note that the theorem is still valid because the additional factor of $e^{-\frac{a^2}{4}t}$  will  cancel when we apply the boundary conditions.
\end{example}
\begin{example}
An important class of stochastic processes are squared Bessel processes, see \cite{FK05}. For the three dimensional case we have
\begin{align}\label{Kolbckeqn}
u_t=2xu_{xx}+3u_x, \ x\geq b>0, t>0.
\end{align}
The elementary solutions  are
\begin{align*}
w_1(x,t;\l)=\frac{1}{\sqrt{x}}e^{-\l^2t}\cos(\l(\sqrt{2x}-\sqrt{2b})),
\end{align*}
and
\begin{align*}
w_2(x,t;\l)=\frac{1}{\sqrt{x}}e^{-\l^2t}\sin(\l(\sqrt{2x}-\sqrt{2b})).
\end{align*}
\end{example}
\begin{example}
We consider the PDE
\begin{align}
u_t=u_{xx}+2(\tanh x)u_x,\ x>0, t>0,
\end{align}
The elementary solutions are
\begin{align*}
w_1(x,t;\xi)&=e^{-\xi^2t}\frac{\cos(\xi x)}{\cosh x}, \
w_2(x,t;\xi) =e^{-\xi^2t}\frac{\sin(\xi x)}{\cosh x}.
\end{align*}
The solution to the BVP will be of the form \eqref{solution}.

\end{example}
We also remark that the PDE $u_t=u_{xx}+2(\coth x) u_x$  arises in hyperbolic geometry. See the discussion in \cite{Dav89}.  The elementary solutions
are
\begin{align*}
w_1(x,t;\xi)&=e^{-\xi^2t}\frac{\cos(\xi x)}{\sinh x}, \
w_2(x,t;\xi) =e^{-\xi^2t}\frac{\sin(\xi x)}{\sinh x}.
\end{align*}
To avoid  singularities one must place the lower boundary at $x=b>0$.

\begin{example}
Next we consider a family of PDEs of the form
\begin{align*}
u_t=u_{xx}+(2x^2+a)u_x+(x^4+ax^2+2x)u,\ x>0, t>0.
\end{align*}
The elementary solutions that we use are
\begin{align*}
w_1(x,t;\xi)&=e^{-\frac14(\xi^2+a^2)t}e^{-\frac13x^3-\frac12ax}\cos\left(\frac{x\xi}{2}\right),\\
w_2(x,t;\xi)&=e^{-\frac14(\xi^2+a^2)t}e^{-\frac13x^3-\frac12ax}\sin\left(\frac{x\xi}{2}\right).
\end{align*}

\end{example}

Variants of the Robin boundary condition can also be studied and we present some examples.

\begin{example}

Naturally we cannot obtain a solution where none exists. Boundary conditions must be suitable to the PDE. We illustrate this with the
equation
\begin{align}
u_t=(x^2-1)u_{xx}+xu_x,\ x>1, t>0.
\end{align}
We require $u(1,t)=f(x),$ $\alpha u(1,t)+\beta\sqrt{x^2-1}u_x(x,t)\big|_{x\to 1}=g(t).$ It is not clear that the standard Robin problem at $x=1$ even has a solution. This is clearly not covered by Theorem \ref{AGentheorem}. Nevertheless the same construction works as in that result.
Our elementary solutions are
\begin{align*}
w_1(x,t;\xi)&=e^{-\xi^2t}\cos\left(\xi\cosh^{-1}{x}\right), \
w_2(x,t;\xi) =e^{-\xi^2t}\sin\left(\xi\cosh^{-1}{x}\right).
\end{align*}
From the initial data we obtain $\varphi(\xi)=F(\xi)+\mathcal{H}_o\psi)(\xi)$ with $F(\xi)=\frac{2}{\pi}\int_0^\infty f(\cosh y)\cos(\xi y)dy$.
The boundary condition gives $$\alpha\varphi(\xi)+\beta\xi\psi(\xi)=2\xi\mathcal{L}^{-1}[g(t);\xi^2],$$ and we have solved equations like these before.
\end{example}

\subsubsection{Extension to time dependent coefficients}
We can also study problems with time dependent coefficients. A full discussion would be lengthy, so we only make some brief remarks.
Suppose that $r:[0,\infty)\to (0,\infty) $ is a continuous positive function. Suppose that the PDE \eqref{genrobthm} has solutions given by \eqref{f3} and \eqref{f4}.
Then the PDE
\begin{align}\label{tdepend}
\frac{1}{r(t)}u_t=P(x)u_{xx}+Q(x)u_x+R(x)u,\ x>b>0,
\end{align} has solutions given by
\begin{align}
\bar{v}_1(x,t;\xi)=\rho(x)\cos(\xi(\mu(x)-\mu(b)))e^{-\xi^2\int_0^t r(s)ds},\label{f5}\\
\bar{v}_2(x,t;\xi)=\rho(x)\sin(\xi(\mu(x)-\mu(b)))e^{-\xi^2\int_0^t r(s)ds}\label{f6}.
\end{align}
The extension of Theorem \ref{genrobthm} to \eqref{tdepend} subject to the boundary condition \eqref{genrobcond} is entirely straightforward and we
leave it the interested reader. See however the remark in Section \ref{finsec}.

Many classes of equations with time dependent coefficients can be solved and we make no attempt to list them all here, but we simply
present as an example equations of the form
\begin{align}\label{tdepend2}
u_t=\sigma(t)u_{xx}+(A(t)-B(t)x)u_x, \ x>0, \sigma(t)>0.
\end{align}
Here $\sigma,A,B$ are continuous. Equations of this type admit solutions of the form
\begin{align}
u_c(x,t)&=\exp(-\xi^2 k(t))\cos\left(\xi(x e^{-\int_0^t B(s)ds}+C(t))\right)\\
u_s(x,t)&=\exp(-\xi^2 k(t))\sin\left(\xi (x e^{-\int_0^t B(s)ds}+C(t))\right),
\end{align}   where $k(t)=\int_0^t\sigma(y)e^{-2\int_0^y B(s)ds}dy$ and $C(t)=\int_0^t A(y)e^{-\int_0^yB(s)ds}dy.$ A variety of boundary value problems for \eqref{tdepend2} can be solved by our method.

We take our solutions to be of the form
\begin{align*}
u(x,t)&=\int_0^\infty \varphi(\xi)\exp(-\xi^2 k(t))\cos\left(\xi (x e^{-\int_0^t B(s)ds}+C(t))\right)d\xi\\&+
\int_0^\infty \psi(\xi)\exp(-\xi^2 k(t))\sin\left(\xi (x e^{-\int_0^t B(s)ds}+C(t))\right)d\xi.
\end{align*}
Since $k(0)=C(0)=0$ we have
\begin{align*}
u(x,0)=\rho(x)(\widehat{\varphi}_c(x)+\widehat{\psi}_s(x))=f(x).
\end{align*}
This is the same condition that we had before. There are various cases where the equations resulting from the boundary conditions can also be solved using the odd and even Hilbert transforms, however some cases are challenging. If $A(t)=0$, then the Neumann and Dirichlet problems are reasonably straightforward. For example, the inhomogeneous Dirichlet problem yields
\begin{align}
u(0,t)=\int_0^\infty \varphi(\xi)e^{-\xi^2 k(t)}d\xi=g(t).
\end{align}
Setting $\xi^2=z$ we have the Laplace transform $\tilde{\Phi}(k(t))=g(t)$. Here $\tilde{\Phi}$ is the Laplace transform of $\tilde{\varphi}(z)=\frac{\varphi(\sqrt{z})}{2\sqrt{z}}.$ Since $k$ is increasing it is invertible and one writes
$\tilde{\Phi}(s)=g\left(k^{-1}(s)\right)$ and this can be inverted explicitly for a large number of cases. This gives us a solution of the Dirichlet problem.

\subsection{A Modified Robin Problem for the Harmonic Oscillator}\label{LT}
The equation for the harmonic oscillator is,
\begin{align}
u_t=\sigma u_{xx}-\mu x^2u,\ \ \ x\in \Omega\subseteq \mathbb{R},\mu>0.
\end{align}
This plays an important role in quantum mechanics, though it would normally be in the form $iu_t=\sigma x_{xx}-\mu x^2u.$ The real form also
plays an important role in the theory of heat kernels. See Davies' book \cite{Dav89}.

It also arises in stochastic analysis. By the Feynman-Kac formula the functional $u(x,t)=\mathbb{E}\left[f(X_t)e^{-\mu\int_0^t X_s^2ds}\big|X_0=x\right]$, where $X=\{X_t;t\geq 0\}$ is a Brownian motion, is given by the solution of the problem
\begin{align*}
u_t&=\frac12 u_{xx}- \mu x^2u, x\in \mathbb{R},\\
u(x,0)&=f(x).
\end{align*}
If $f(x)=e^{-\l x}$, then the solution of the initial value problem gives the Laplace transform of the joint density of $(X_t,\int_0^t X_s^2ds)$.
Our intention is to solve a modified Robin problem which has a special form. However we note here  that the Dirichlet and Neumann problems   can easily be solved by our method. In fact there are numerous problems problems that we can  solve.  The standard Robin problem reduces to the solution of a different type of integral equation of Laplace transform convolution type. We will present it elsewhere.

We will solve the problem
\begin{align*}
u_t&=  u_{xx}- x^2u,\ x\ \geq 0,\\
u(x,0)&=f(x),\\
\alpha u(0,t)+\beta e^{-2t}u_x(0,t)&=g(t).
\end{align*}
We will suppose that $\alpha$ and $\beta$ are non zero.  It is worth noting that if we map  this to the heat equation, we obtain a very complicated moving boundary problem which appears to be intractable.

The boundary condition describes a situation where the lower boundary starts off as partially absorbing, partially reflecting  and exponentially decays to a purely absorbing boundary. Such a situation can arise in many settings, such as the design of materials which reflect, say, alpha or beta particles. These can cause serious damage to living tissue.

We can let $u$ be the amount of radioactive material that has been absorbed.  As the protective material breaks down the proportion reflected decreases to zero and so the boundary will become purely absorbing. For small $t$, the reflectivity will decay linearly, but after a certain point is reached, the reflectivity will evaporate exponentially fast. These are important considerations in the design of systems to protect against radiation.
There is a very large literature on this subject and we can only suggest a survey such as \cite{But77}. Our problem may be regarded as a toy problem that could be of value. However, we present it here purely for its mathematical interest. We prove the following theorem.
\begin{thm}
Let $f\in L^1([0,\infty)$, $g(t)=\int_0^\infty G(x)e^{-xt}d x$, where  $G\in L^1\left([0,\infty),\frac{d x }{\Gamma\left(\frac{1+x}{4}\right)}\right)$. Then the boundary value problem
\begin{align*}
u_t&=  u_{xx}- x^2u,\ x\ \geq 0,\\
u(x,0)&=f(x),\\
\alpha u(0,t)+\beta e^{-2t}u_x(0,t)&=g(t),
\end{align*}
where $\alpha,\beta$ are nonzero, has a solution given by
\begin{align}
u(x,t)&=\int_0^\infty\varphi(\l)w_1(x,t;\l)d \l+\int_0^\infty\psi(\l)w_2(x,t;\l)d \l.
\end{align}
Here
\begin{align*}
w_1(x,t;\l)&=\exp\left(t+\frac12x^2-\frac14\l^2(e^{4t}-1)\right)\cos(\l xe^{2t}),\\
w_2(x,t;\l)&=\exp\left(t+\frac12x^2-\frac14\l^2(e^{4t}-1)\right)\sin(\l xe^{2t}),
\end{align*}
\begin{align*}
\psi(\l)=\frac{\alpha\left[\alpha(\mathcal{H}_eF)(\l)-(\mathcal{H}_eK)(\l)+\frac{\beta\l  }{\alpha}K(\l)-\l\beta F(\l)\right]}{\alpha^2+\beta^2\l^2},
\end{align*}
$\varphi(\l)=\frac{1}{\alpha}K(\l)-\frac{\l\beta}{\alpha}\psi(\l),$
$F(\l)=\frac{2}{\pi}\int_0^\infty e^{-\frac12y^2}f(y)\cos(\l y)d y$ ,
$K(\l)=\frac{\l}{2}\mathcal{L}^{-1}\left[\bar{g}(s);\frac{\l^2}{4}\right]$ and
$\bar{g}(s)=\frac{1}{(1+s)^{\frac14}}g\left(\frac14\ln(1+s)\right).$
\end{thm}
\begin{proof}
We begin with the following solutions of the PDE
\begin{align}
w_1(x,t;\l)&=\exp\left(t+\frac12x^2-\frac14\l^2(e^{4t}-1)\right)\cos(\l xe^{2t}),\\
w_2(x,t;\l)&=\exp\left(t+\frac12x^2-\frac14\l^2(e^{4t}-1)\right)\sin(\l xe^{2t}).
\end{align}
It is straightforward to check that these satisfy the equation for every $\l$. (Notice that they are not however of the form \eqref{f1} or \eqref{f2}).

As usual we construct a new solution of the PDE by setting
\begin{align}
u(x,t)&=\int_0^\infty\varphi(\l)w_1(x,t;\l)d \l+\int_0^\infty\psi(\l)w_2(x,t;\l)d \l.
\end{align}
The initial condition then implies that
\begin{align}
e^{\frac12x^2}\left[\widehat{\varphi}_c(x)+\widehat{\psi}_s(x)\right]&=f(x).
\end{align}
Reasoning as previously, we deduce that
\begin{align}
\varphi(\l)=\frac{2}{\pi}\int_0^\infty e^{-\frac12y^2}f(y)\cos(\l y)d y+(\mathcal{H}_o\psi)(\l).
\end{align}
We set $F(\l)=\frac{2}{\pi}\int_0^\infty e^{-\frac12y^2}f(y)\cos(\l y)d y.$

Now we easily see that
\begin{align}
u(0,t)&=\int_0^\infty\varphi(\l)\exp\left(t-\frac14\l^2(e^{4t}-1)\right)d \l
\end{align}
and
\begin{align}
u_x(0,t)&=\int_0^\infty\l\psi(\l)\exp\left(3t-\frac14\l^2(e^{4t}-1)\right)d \l.
\end{align}

From the boundary condition we have
\begin{align}\label{bceqn1}
\alpha \int_0^\infty\varphi(\l)e^{ -\frac{\l^2} {4}(e^{4t}-1) }d \l+\beta \int_0^\infty\l\psi(\l)e^{-\frac{\l^2} {4}(e^{4t}-1)}d \l=e^{-t}g(t).
\end{align}

Now we set $s=e^{4t}-1$, which implies $t=\frac14\ln(1+s).$ We also put $\l^2=4z.$ This transforms  \eqref{bceqn1} to
\begin{align}
\alpha \int_0^\infty\frac{\varphi(2\sqrt{z})}{\sqrt{z}}\exp\left(-zs\right)d z+2\beta \int_0^\infty \psi(2\sqrt{z})\exp\left(-zs\right)d z=\bar{g}(s)
\end{align}
in which $$\bar{g}(s)=\frac{1}{(1+s)^{\frac14}}g\left(\frac14\ln(1+s)\right).$$

Inverting the Laplace transform gives
\begin{align}
\alpha \frac{\varphi(2\sqrt{z})}{\sqrt{z}} +2\beta  \psi(2\sqrt{z}) =\mathcal{L}^{-1}\left[\bar{g}(s);z\right],
\end{align}
and so
\begin{align}
\alpha\varphi(\l)+\beta\l\psi(\l)=\frac{\l}{2}\mathcal{L}^{-1}\left[\bar{g}(s);\frac{\l^2}{4}\right]=K(\l),
\end{align}
which immediately yields
\begin{align}
\alpha(F(\l)+(\mathcal{H}_o\psi)(\l))=K(\l)-\beta\l\psi(\l).
\end{align}

Taking the even Hilbert transform of both sides produces the result
\begin{align}
\alpha(\mathcal{H}_eF)(\l)-\alpha\psi(\l))=(\mathcal{H}_eK)(\l)-\beta\l(\mathcal{H}_o\psi)(\l).
\end{align}

However we also have

$$(\mathcal{H}_o\psi)(\l)=\frac{1}{\alpha}\left[K(\l)-\beta\l\psi(\l)\right]-F(\l).$$
Whence

\begin{align}
\alpha(\mathcal{H}_eF)(\l)-\alpha\psi(\l))=(\mathcal{H}_eK)(\l)-\beta\l\left[\frac{1}{\alpha}\left[K(\l)-\beta\l\psi(\l)\right]-F(\l)\right].
\end{align}
Solving for $\psi$ gives

\begin{align}
\psi(\l)=\frac{\alpha\left[\alpha(\mathcal{H}_eF)(\l)-(\mathcal{H}_eK)(\l)+\frac{\beta\l  }{\alpha}K(\l)-\l\beta F(\l)\right]}{\alpha^2+\beta^2\l^2}.
\end{align}
From this we can find $\varphi.$ Proving convergence of the integrals for these choices of $\varphi$ and $\psi$ proceeds much as in our previous example. Since the details are similar we focus on the function $K$ which has a different structure than in the previous case. To proceed we suppose that
\begin{align}
g(s)=\int_0^\infty G(x)e^{-xs}d x.
\end{align}
Then we immediately see that
\begin{align}
\frac{1}{(1+s)^{\frac14}}g\left(\frac14\ln(1+s)\right)=\int_0^\infty \frac{G(x)}{(1+s)^{\frac{1+x}{4}}}d x.
\end{align}

Taking the inverse Laplace transform we obtain
\begin{align}
\mathcal{L}^{-1}\left[\frac{1}{(1+s)^{\frac14}}g\left(\frac14\ln(1+s)\right);z\right]=\int_0^\infty\frac{  G(x)}{\Gamma\left(\frac{1+x}{4}\right)}e^{-\frac{(x+1)z}{4}}
z^{\frac{x}{4}-\frac{3}{4}}d x.
\end{align}
See \cite{RK66} for the inverse Laplace transform.

Thus we obtain the expression
\begin{align}\label{intK}
K(\l)&=\left(\frac{\l}{2}\right)^{-\frac{1}{2}}e^{-\frac{\l^2}{4}}\int_0^\infty\frac{  G(x)}{\Gamma\left(\frac{1+x}{4}\right)}e^{-\frac{x\l^2}{16}}
\left(\frac{\l}{2}\right)^{\frac{x}{2}}d x.
\end{align}

If we suppose that $G\in L^1\left([0,\infty),\frac{d x }{\Gamma\left(\frac{1+x}{4}\right)}\right)$, then an application of the Dominated Convergence Theorem gives
\begin{align}
\lim_{\l\to \infty}\int_0^\infty\frac{1}{\Gamma\left(\frac{1+x}{4}\right)}G(x)e^{-\frac{x\l^2}{16}}
\left(\frac{\l}{2}\right)^{\frac{x}{2}}d x=0.
\end{align}
We also have
\begin{align}
\int_0^\infty\frac{1}{\Gamma\left(\frac{1+x}{4}\right)}G(x)e^{-\frac{x\l^2}{16}}
\left(\frac{\l}{2}\right)^{\frac{x}{2}}d x\big|_{\l=0}=0.
\end{align}
Now the Lebesgue integral of a measurable function returns a uniformly continuous function. See \cite{SS2005}. Hence the integral in (\ref{intK}) is bounded and   we have the inequality
\begin{align}
|K(\l)|\leq \frac{C}{\sqrt{\l}}e^{-\frac{\l^2}{4}},
\end{align}
for some positive constant $C$ depending on $G$. Thus the even Hilbert transform of $K$ exists. We can prove that $(\mathcal{H}_eK)(\l)w_{k}(x,t;\l)$
is integrable for $k=1,2$, by similar means to the Black-Scholes equation case. If we suppose that $f$ is integrable then similar statements can be made for the even Hilbert transform of $F$. The remainder of the proof is similar to our previous examples.
\end{proof}

\section{The Robin Problem for a Five Dimensional Squared Bessel Process}\label{BESQ5}
The method is effective for elementary solutions of greater complexity. For example solutions of the form
\begin{align}
w_1(x,t)&=(\rho_1(x,t;\xi)\cos(\xi A_1(x,t))+\rho_2(x,t;\xi)\sin(\xi A_2(x,t)))e^{-\xi^2t},\label{ze1}\\
w_2(x,t)&=(\rho_3(x,t;\xi)\cos(\xi A_1(x,t))+\rho_4(x,t;\xi)\sin(\xi A_2(x,t)))e^{-\xi^2t}\label{ze2}
\end{align}
can be used. It is possible to formulate an analogue of Theorem \ref{AGentheorem} for equations with elementary solutions of this form.
However we omit this for brevity. Instead we will content ourselves by working through the details of a particular example of interest.
We solve the problem

\begin{align}
u_t&=2xu_{xx}+5u_x,\ x>b, t>0,\label{5dsbes}\\
u(x,0)&=f(x), x>b,\\
\alpha u(b,t)+\beta u_x(b,t)&=g(t).
\end{align}
We assume that $f$ is continuously differentiable. This is associated with a five dimensional squared Bessel process and is not covered by Theorem \ref{AGentheorem}. However our method can be extended to
cover it.
We look for a solution which can be written
\begin{align}
u(x,t)&=\int_0^\infty  \varphi (\xi )w_1(x,t;\xi)d\xi + \int_0^\infty \psi (\xi )w_2(x,t;\xi)d\xi,
\end{align}
in which
\begin{align}
w_1(x,t;\xi)&=e^{-\xi^2t} x^{-3/2} \left(2 \xi  \sqrt{x} \sin\left(\xi m(x)\right)+\sqrt{2} \cos \left(\xi m(x)\right)\right) , \label{s1}\\
w_2(x,t;\xi)&=e^{-\xi^2t}x^{-3/2} \left(\sqrt{2} \sin\left(\xi m(x)\right)-2 \xi  \sqrt{x} \cos \left(\xi m(x)\right)\right) .\label{s2}
\end{align}
Here $m(x)=\sqrt{2x}-\sqrt{2b}.$

If we make the change of variables $x=\frac12(y+\sqrt{2b})^2$ the initial condition reduces to
\begin{align*}
&\int_0^\infty\varphi(\xi)(\cos(\xi y)+(y+\sqrt{2b})\xi\sin(\xi y))d\xi\\&\quad\quad+\int_0^\infty\psi(\xi)(\sin(\xi y)-(y+\sqrt{2b})\xi\cos(\xi y))d\xi\\
&=\frac{1}{\sqrt{2}}\left(\frac12(y+\sqrt{2b})^2\right)^{3/2}f\left(\frac12(y+\sqrt{2b})^2\right) =\tilde{f}(y).
\end{align*}
This seems to be quite different to the previous cases that we have encountered. However we can  rewrite the above as
\begin{align}
\left(1-(y+\sqrt{2b})\frac{d}{dy}\right)(\widehat{\varphi}_c(y)+\widehat{\psi}_s(y))=\tilde{f}(y).
\end{align}
Or equivalently
\begin{align}
(\frac{d}{dy}-\frac{1}{y+\sqrt{2b}})(\widehat{\varphi}_c(y)+\widehat{\psi}_s(y))&=-\frac{1}{y+\sqrt{2b}}\tilde{f}(y).
\end{align}

Introducing the integrating factor $\displaystyle{\frac{1}{y+\sqrt{2b}}}$ we have
\begin{align}
\frac{d}{dy}\left(\frac{1}{y+\sqrt{2b}}\left( \widehat{\varphi}_c(y)+\widehat{\psi}_s(y) \right)\right)=-\frac{1}{(y+\sqrt{2b})^2}\tilde{f}(y).
\end{align}
Hence
\begin{align}\label{eff}
\widehat{\varphi}_c(y)+\widehat{\psi}_s(y)= (y+\sqrt{2b})\left(\int_0^y\frac{-1}{(\eta+\sqrt{2b})^2}\tilde{f}(\eta)d\eta+I\right)=F(y),
\end{align}
where $I$ is a constant of integration. To determine the value of $I$ we need an initial condition and we notice that
$\widehat{\psi}_s(0)=0$ and $\widehat{\varphi}_c(0)=\int_0^\infty\varphi(\xi)d\xi.$
Thus we need the value of the integral $\int_0^\infty\varphi(\xi)d\xi.$ This same constant also arises from the boundary conditions. We will see that it can be obtained by solving a pair of simultaneous equations.

From \eqref{eff} we obtain the familiar condition
\begin{align}\label{famcond}
\varphi(\xi)=\widehat{F}(\xi)+(\mathcal{H}_o\psi)(\xi),
\end{align}
with $\widehat{F}(\xi)=\frac{2}{\pi}\int_0^\infty F(\eta)\cos(\eta\xi)d\eta.$

We remark here that for the squared
Bessel processes of odd order higher than five,   the  condition \eqref{famcond} is obtained by solving equations of Euler type for $\widehat{\varphi}_c(y)+\widehat{\psi}_s(y)$
and the initial conditions are found by solving a system of linear equations. However the details are quite involved and so will be presented elsewhere.

The Robin boundary condition gives us
\begin{align}\label{BCSQ5}
\frac{ ((2b\alpha-3\beta)+2b\beta\xi^2)\varphi(\xi)}{2b^{5/2}}+\frac{3\beta-2b\alpha}{b^2}\xi\psi(\xi)=2\xi\mathcal{L}^{-1}[g(t);\xi^2]=G(\xi).
\end{align}
This is equivalent to
\begin{align}\label{fon}
(A+B\xi^2)(\mathcal{H}_o\psi)(\xi) +C\xi\psi=G(\xi)-(A+B\xi^2)\widehat{F}(\xi),
\end{align}
with $A=\frac{(2b\alpha-3\beta)}{2b^{5/2}}$, $C=\frac{3\beta-2b\alpha}{b^2}$. We assume that $A+B\xi^2$ has no real roots.
We take the even Hilbert transform of both sides, just as before. This gives
\begin{align}\label{gilp}
-A\psi+B\left(-\frac{2\xi}{\pi}P\int_0^\infty\mathcal{H}_o\psi(\eta)d\eta-\xi^2\psi(\xi)\right)+ C\xi\left(\mathcal{H}_o\psi\right)(\xi)=J(\xi)
\end{align}
where $J(\xi)=(\mathcal{H}_eK)(\xi)$ and $K(\xi)= G(\xi)-(A+B\xi^2)\widehat{F}(\xi) .$ Now as before, we require
\begin{align}
P\int_0^\infty\left(\mathcal{H}_o\psi\right)(\eta)d\eta&=P\int_0^\infty(\varphi(\eta)-\widehat{F}(\eta))d\eta.
\end{align}
We will drop the principal value by assuming suitable integrability. The procedure for determining the value of the integral $\int_0^\infty \varphi(\xi)d\xi$ will be addressed  below.

From \eqref{fon} we see that
\begin{align}
(\mathcal{H}_o\psi)(\xi)&=\frac{K(\xi)-C\xi\psi(\xi)}{A+B\xi^2}.
\end{align}
Together with \eqref{gilp} we have a pair of simultaneous equations for $\psi$ and $\mathcal{H}_o\psi $ and we obtain
\begin{align}\label{gilp2}
-\frac{(A+B\xi^2)^2+C^2\xi^2}{A+B\xi^2}\psi(\xi) =(\mathcal{H}_eK)(\xi)+B\frac{2\xi}{\pi} \int_0^\infty\mathcal{H}_o\psi(\eta)d\eta-k(\xi),
\end{align}
where $k(\xi)=\frac{C\xi K(\xi) }{A+B\xi^2}.$ This gives
\begin{align}
\psi(\xi)=-(A+B\xi^2)\frac{(\mathcal{H}_eK)(\xi)+B\frac{2\xi}{\pi}\left( \int_0^\infty(\varphi(\eta)-\widehat{F}(\eta))d\eta\right)-\frac{C\xi K(\xi) }{A+B\xi^2}}{(A+B\xi^2)^2+C^2\xi^2}.
\end{align}
From this we obtain $\varphi.$ To complete the calculation we need to determine $\int_0^\infty\varphi(\xi)d\xi.$
First from the construction of the solution we have
\begin{align}
u(b,0)&=f(b)=b^{-3/2}\sqrt{2}\int_0^\infty\varphi(\xi)d\xi-\frac{2}{b}\int_0^\infty\xi\psi(\xi)d\xi.
\end{align}
Integrating both sides of \eqref{BCSQ5} gives
\begin{align}
A\int_0^\infty\varphi(\xi)d\xi+B\int_0^\infty\xi^2\varphi(\xi)d\xi+C\int_0^\infty\xi\psi(\xi)d\xi=\int_0^\infty G(\xi)d\xi.
\end{align}
We next calculate $u_x(b,0)$ and after rearranging obtain
\begin{align}
\int_0^\infty\xi\psi(\xi)d\xi&=-\frac{\sqrt{2a}}{3}\int_0^\infty \xi^2\varphi(\xi)d\xi+\frac{1}{\sqrt{2b}}\int_0^\infty\varphi(\xi)d\xi+f'(b).
\end{align}
Basic algebra leads to the pair of simultaneous equations
\begin{align*}
H_1\int_0^\infty\varphi(\xi)d\xi+H_2\int_0^\infty\xi^2\varphi(\xi)d\xi&=\frac{b^{3/2}}{\sqrt{2}}f(b)
-\frac{\sqrt{2b}}{C}\int_0^\infty G(\xi)d\xi,\\
H_3\int_0^\infty\varphi(\xi)d\xi+H_4\int_0^\infty\xi^2\varphi(\xi)d\xi&=
\frac{b^2}{3}f'(b)-\int_0^\infty\frac{G(\xi)}{C}d\xi.
\end{align*}
 $H_1=\left(1+\frac{\sqrt{2b} A}{C}\right)$, $H_2=\frac{\sqrt{2b}B}{C}$, $H_3=\left(\sqrt{\frac{b}{2}}+\frac{A}{C}\right)$, $H_4=\left(\frac{B}{C}-\frac{\sqrt{2b}}{3}\right)$.
This gives $\int_0^\infty\varphi(\xi)d\xi=\gamma$, where
\begin{align}
\gamma&=\frac{H_4\left(\frac{b^{3/2}}{\sqrt{2}}f(b)
-\frac{\sqrt{2b}}{C}\int_0^\infty G(\xi)d\xi \right)-H_2\left(\frac{b^2}{3}f'(b)-\int_0^\infty\frac{G(\xi)}{C}d\xi\right) }{H_1H_4-H_2H_3}.
\end{align}
Where we assume that $H_1H_4-H_2H_3\neq 0.$ If $g(t)=0$ then this simplifies to
\begin{align}
\int_0^\infty\varphi(\xi)d\xi=\frac{\frac{b^{3/2}}{\sqrt{2}}H_4f(b)
 -\frac{b^2}{3}H_2f'(b)  }{ H_1H_4-H_2H_3}.
\end{align}
We immediately obtain $$F(y)=(y+\sqrt{2b})\left(\int_0^y\frac{-1}{(\eta+\sqrt{2b})^2}\tilde{f}(\eta)d\eta+\frac{\gamma}{\sqrt{2b}}\right).$$

The question then arises as to how we are to interpret the Fourier cosine transform of $y+\sqrt{2b}?$ This is in terms of distributions.
Clearly $\int_0^\infty\delta(\xi)\cos(\xi y)d\xi=1.$ Hence $\frac{2}{\pi}\int_0^\infty\cos(\xi y)dy=\delta(\xi).$  Now let $\phi$ be a suitable
test function of Schwartz class. We  define the distribution
$\Lambda(\phi)=\int_0^\infty y\widehat{\phi}_c(y)dy$ and as a distribution $\left(\int_0^\infty y \cos(\xi)dy,\phi\right)=\Lambda(\phi).$
Note that the solutions \eqref{s1} and \eqref{s2} are in the Schwartz class, so that $\Lambda(w_1)$ and $\Lambda(w_2)$ are well defined
and can be computed easily.

It is possible to give conditions which guarantee convergence and to prove that we do indeed have a solution
of the Robin problem. This proceeds along the lines of the treatment given in the cases of the Black-Scholes equation and the harmonic oscillator. To avoid repeating the same arguments as before, we will omit this analysis. However we are able to construct a solution to the problem by this method, which again does not require us to know the fundamental solution.


\section{Future Directions}\label{finsec}

This paper has introduced a new method which yields new representations for the solution of important boundary value problems. The idea of constructing solutions to boundary value problems from elementary solutions of the PDE, without use of a fundamental solution is potentially very important. Naturally there are many open problems.

First we would like a fuller characterisation of the types of equations for which the method is effective. We have one result  along those lines, and others are possible. We can give a characterisation of PDEs which possess solutions of the form \eqref{ze1} and \eqref{ze2}. We proceed as in the discussion at the start of Section \ref{QE2d}. One substitutes the solutions into an arbitrary linear PDE and this leads to conditions which guarantee that we have a solution.

There are also special cases that we have not gone into in detail, in order to keep the current paper to a manageable length. For most of these cases, the analysis is basically the same as that presented for our main problems. In many cases they are easier because certain terms disappear from the equations. For example in the derivation of Theorem 3.3, assuming that $\left(\alpha+\frac{\mu\beta}{b}\right)=0$  significantly reduces the difficulty in determining
$\varphi$ and $\psi.$ One simply inverts a single Laplace transform to obtain $\psi$ and $\varphi$ can be determined from this.

We do not yet know the full range of equations with time dependent coefficients which can be studied by this technique. However we made some preliminary remarks on this earlier. If we consider equation \eqref{tdepend}, the extension of Theorem \ref{AGentheorem} to equations of this form is not at all difficult. We simply construct a solution
\begin{align*}
u(x,t)&=\int_0^\infty\varphi(\xi)\rho(x)\cos(\xi(\mu(x)-\mu(b)))e^{-\xi^2\int_0^t r(s)ds}d\xi\\&\quad  +\int_0^\infty\psi(\xi)\rho(x)\sin(\xi(\mu(x)-\mu(b)))e^{-\xi^2\int_0^t r(s)ds}d\xi.
\end{align*}
Note that $k(t)=\int_0^t r(s)ds$ is increasing and so is invertible. One has
\begin{align*}
u(x,0)=\rho(x)\left(\widehat{\varphi}_c(\mu(x)-\mu(b))+\widehat{\psi}_s(\mu(x)-\mu(b))\right)=f(x),
\end{align*}
and the equations arising from the boundary conditions can be solved by the use of the odd and even Hilbert transforms.
The analysis is similar to the case when $r(s)=1$ and the solution is a modification of that given in Theorem \ref{AGentheorem}.
The solution is in terms of the inverse Laplace transform of $g(k^{-1}(s))$, as with our example for the equation \eqref{tdepend2}.

For \eqref{tdepend2}, it is clear that there are numerous possible cases that can be considered. We have seen that the Dirichlet problem is straightforward if one sets $A(t)=0$, which corresponds to a time dependent Ornstein-Uhlenbeck process. The Neumann problem can also be solved in the same manner. For different choices of the coefficients, a variety of boundary value problems can be solved. There are other time dependent equations that can be studied, but we will not discuss them here.

Another application of these techniques is to equations with nonstandard boundary conditions. We mentioned previously that certain moving boundary problems can be solved. We will give one example. Consider the PDE $u_t=u_{xx}-xu$, $x\geq t^2$, subject to $u(t^2,t)=g(t),$ $u(x,0)=f(x).$ The PDE has solutions
\begin{align}
u_1(x,t;\xi)=e^{\frac13t^3-xt-\xi^2t}\cos(\xi(x-t^2)),\\
u_2(x,t;\xi)=e^{\frac13t^3-xt-\xi^2t}\sin(\xi(x-t^2)).
\end{align}
We can solve a number of different problems with these elementary solutions, see \cite{Cra2008}. To solve the moving boundary problem
we define
\begin{align}
u(x,t)=\int_0^\infty\varphi(\xi)u_1(x,t;\xi)d\xi+\int_0^\infty\varphi(\xi)u_2(x,t;\xi)d\xi.
\end{align}
Then $u(x,0)=\widehat{\varphi}_c(x)+\widehat{\psi}_s(x)=f(x).$ The moving boundary condition leads to
\begin{align}
\int_0^\infty\varphi(\xi)e^{-\xi^2t}d\xi=e^{2t^3/3}g(t).
\end{align}
This reduces to a Laplace transform and inverting the transform gives us $\varphi$ and from this we obtain $\psi.$ There are other moving boundary problems that one can study for different PDEs. This is ongoing work.

We also remark that different families of elementary solutions which do not involve sines and cosines can also be used in our basic construction. There are also PDEs which possess elementary solutions for which the integral equations arising from \eqref{formsoln} can be solved by means other than the Hilbert transform. Much work remains to be done on these problems, however see the preprint \cite{Crad2019A} for work in this direction.
\section{Appendix 1}
The classical approach to solving the Robin problem when $g$ is nonzero is presented here. For a more technical treatment of the solution of boundary problems for parabolic operators, we recommend \cite{Fried64}. We wish to solve
\begin{align}
 u_t(S,t)& = L u(S,t),\label{m1} \\
u(S,0) &= f(S), \label{m2}  \\
\alpha u(b,t) + \beta u_S(b,t) + \gamma u_{SS}(b,t)&= g(t),\label{m3}
   \end{align}
with $S\geq b, t>0, b>0$ and where $L$ is a positive second order operator in $S$. In our case, $L$ is the second order Black-Scholes operator, namely
\begin{equation*}
    Lu(S,t) = \frac{1}{2} \sigma^2 S^2 u_{SS}(S,t) + r S u_S(S,t).
\end{equation*}
In order to solve problem \eqref{m1}-\eqref{m3}, we set $u(S,t)=v(S,t)+h(S,t)$, where
\begin{equation}
    \alpha v(b,t) + \beta v_S(b,t) + \gamma v_{SS}(b,t) = 0,
\nonumber\end{equation}
and
\begin{equation}
\label{h function req}
    \alpha h(b,t) + \beta h_S(b,t) + \gamma h_{SS}(b,t) = g(t).
\end{equation}
This gives
\begin{equation}
   \begin{cases}
        v_t(S,t) = L v(S,t) + K(S,t) \\
        v(S,0) = f(S) - h(S,0)\\
        \alpha v(b,t) + \beta v_S(b,t) + \gamma v_{SS}(b,t) = g(t),
   \end{cases}
\end{equation}
where $K(S,t)=Lh-h_t(S,t)$ and $S>0, t>0$. We assume $g(t)$ to be differentiable and integrable.   We then have the following result. We omit the proof.
It is available on request.
\begin{prop}
\label{Proposition: representation second order operator}
    Let $h(S,t)$ be a function that satisfies \eqref{h function req} and let $q(S,y,t)$ a fundamental solution of $q_t=Lq$, such that $\alpha q(b,y,t) + \beta q_S(b,y,t) + \gamma q_{SS}(b,y,t) = 0$.
    Then the solution of the
boundary value problem \eqref{m1}-\eqref{m3}, can be written as
\begin{equation}
    \begin{aligned}
        u(S, t)= & h(S, t)+\int_{b}^{\infty}(f(y)-h(y, 0)) q(S, y, t) d y \\
        &+\int_{0}^{t} \int_{b}^{\infty}\left(L h(y, \tau)-\frac{\partial h}{\partial t}(y, \tau)\right) q(S, y, t-\tau) d y d \tau.\nonumber
    \end{aligned}
\end{equation}
\end{prop}
In general, the function $h(S,t)$ is not unique and different choices for it may lead to different representations of the solution to the BVP. The fundamental solution in this case is given by \eqref{trans dens up to delta} in which
\begin{align}
    \begin{split}
        &\tilde{p}(S,y,t)=\frac{S^{\mu } b^{-\frac{\ln \left(\frac{y}{b}\right)+2 \ln (S)}{\sigma ^2 t}-\mu }M(S,t)N(S,t)  }{4 \sqrt{2 \pi \Delta } y \tilde{\beta } \left(\sigma ^2 t\right)^{3/2} \left(\tilde{\alpha }+\mu  \left(\tilde{\beta }+\gamma  \mu \right)\right)} e^{-\frac{y \left(\tilde{\beta }+\sqrt{\Delta
   }\right)}{2 \gamma }}\times \\
& \exp\left(-\frac{4 \left(-4 \ln (b) \ln (S)+\ln ^2(b)+\ln ^2(S)\right)+4 \ln ^2\left(\frac{y}{b}\right)+\left(\sigma ^2-2 r\right)^2 t^2 }{8 \sigma ^2
   t}\right),\nonumber
    \end{split}
\end{align}
where
\begin{align*}
M(S,t)&=\left(-\sqrt{\Delta } \tilde{\beta }+R
   \left(2 \gamma  f(b) \left(\tilde{\alpha }+\mu  \left(\tilde{\beta }+\gamma  \mu \right)\right)-\tilde{\beta } \left(\tilde{\beta }+2 \gamma  \mu \right)\right)\right. \\
   & \left. \tilde{\beta } \left(\tilde{\beta }+2 \gamma  \mu \right) +2
   \sqrt{\Delta } \tilde{\beta } e^{\frac{y \left(\tilde{\beta }+2 \gamma  \mu +\sqrt{\Delta }\right)}{2 \gamma }}+\sqrt{\Delta } \tilde{\beta } \left(-e^{\frac{\sqrt{\Delta }
   y}{\gamma }}\right)\right),
\end{align*}
$R=\left(e^{\frac{\sqrt{\Delta } y}{\gamma }}-1\right) $,
\begin{align*}
N(s,t)&= \left(\sqrt{2 \pi } \sigma ^3 t^{3/2} \tilde{\alpha } b^{\frac{\log \left(\frac{y}{b}\right)}{\sigma ^2 t}} e^{\frac{\ln
   ^2\left(\frac{y}{b}\right)+\ln ^2(b)+\ln ^2(S)}{2 \sigma ^2 t}}Z(S,y)+2 b^{\frac{\ln (S)}{\sigma ^2 t}} S^{-\frac{\ln
   \left(\frac{y}{b}\right)}{\sigma ^2 t}} \times\right. \\
   &\left.\left(b^{\frac{2 \ln \left(\frac{y}{b}\right)}{\sigma ^2 t}} \left(\sigma ^2 t \tilde{\beta }+\gamma  \ln
   \left(\frac{y}{b}\right)\right)+S^{\frac{2 \ln \left(\frac{y}{b}\right)}{\sigma ^2 t}} \left(\sigma ^2 t \tilde{\beta }-\gamma  \ln \left(\frac{y}{b}\right)\right)\right) \right. \\
   & \left.-2
   \gamma  \ln \left(\frac{S}{b}\right) S^{-\frac{\ln \left(\frac{y}{b}\right)}{\sigma ^2 t}} b^{\frac{2 \ln \left(\frac{y}{b}\right)+\ln (S)}{\sigma ^2 t}} -2 \gamma  \ln
   \left(\frac{S}{b}\right) b^{\frac{\ln (S)}{\sigma ^2 t}} S^{\frac{\ln \left(\frac{y}{b}\right)}{\sigma ^2 t}}\right)
\end{align*}
and $Z(S,y)=\text{erf}\left(\frac{\ln \left(\frac{S}{b}\right)+\ln \left(\frac{y}{b}\right)}{\sqrt{2} \sigma
   \sqrt{t}}\right)-\text{erf}\left(\frac{\ln \left(\frac{S}{y}\right)}{\sqrt{2} \sigma  \sqrt{t}}\right)$.






